\documentclass{amsart}

\usepackage{amsfonts}
\usepackage{amsmath}
\usepackage{amssymb}
\usepackage{amsthm}
\usepackage{newlfont}

\newtheorem{lemma}{Lemma}[section]
\newtheorem{pro}{Proposition}[section]
\newtheorem{theorem}{Theorem}[section]
\newtheorem{cor}{Corollary}[section]

\newtheorem{rmk}{Remark}[section]
\newtheorem{exm}{Example}[section]
\allowdisplaybreaks

\begin{document}

\title{Small-time Sampling Behaviour 
of a Fleming-Viot Process}


\author{Youzhou Zhou}
\address{The School of Statistics and Mathematics\\
Zhongnan University of Economics and Law\\ 
182 South Lake Avenue, East Lake New Technology Development Zone, Wuhan, China. 430073}
\email{youzhouzhou2014@znufe.edu.cn}



\date{\today}

\dedicatory{}

\begin{abstract}
The Fleming-Viot process with parent-independent mutation process is one particular neutral population genetic model. As time goes by, some initial species are replaced by mutated ones gradually. Once the population mutation rate is high, mutated species will elbow out all the initial species very quickly. Small time behaviour in this case seems to be the key to understand this fast transition. The small-time asymptotic results related to time scale $\frac{t}{\theta}$ and $a(\theta)t$, where $\lim_{\theta\to\infty}\theta a(\theta)=0$, are obtained in \cite{MR1815182},\cite{MR1649005}, \cite{MR1887170} and \cite{MR2184086}, respectively. Only the behaviour under the scale $t(\theta)$, where $\lim_{\theta\to\infty}\theta a(\theta)=\infty$,
was left untouched. In this paper, the weak limits under various small time scales are obtained. Of particular interest is the large deviations for the small-time transient sampling distributions, which reveal interesting phase transition. Interestingly, such a phase transition is uniquely determined by some species diversity indices. 
\end{abstract}

\maketitle

\section{Introduction}

Neutral population genetic models are concerned with evolution of species under selectively neutral mutation pressure. Let $S$ be the set of all allele types. For simplicity, the type space $S$ is usually taken to be a compact metric space. Let $r(\cdot,\cdot)$ be the metric compatible with the compact topology in $S$. We, therefore, can use a probability measure on $S$ to describe the distribution of allele frequencies.  Denote $\mathcal{P}(S)$ to be the space of Borel probability measures on $S$ and let $C(S)$ be a continuous function space on $S$. We define $\langle f,\mu\rangle=\int_{S} fd\mu$ for $\mu\in\mathcal{P}(S),f\in C(S)$. Then $\mathcal{P}(S)$ is also compact under weak topology generated by metric
 $$
 d_{w}(\mu,\nu)=\sum_{i=1}^{\infty}\frac{|\langle f_{i},\mu\rangle-\langle f_{i},\nu\rangle|\wedge1}{2^{i}}, \mu,\nu\in\mathcal{P}(S), 
 $$
where $\{f_{i},i\geq1\}$ is a countable dense subset in $C(S)$.

   One particular neutral population genetic model is the Fleming-Viot process with parent-independent mutation in \cite{MR1205982}. It is a $\mathcal{P}(S)$-valued diffusion process $Z_{t}$, characterized by generator
\begin{align*}
(\mathcal{L}\varphi)(\mu)=\frac{1}{2}\int_{S}\int_{S}\mu(dx)(\delta_{x}(dy)-\mu(dy))\frac{\delta^{2}\varphi(\mu)}{\delta\mu(x)\delta\mu(y)}+\int_{S}\mu(dx)A\left(\frac{\delta\varphi(\mu)}{\delta\mu(\cdot)}\right)(x),
\end{align*}
where, for a given diffuse probability measure $\nu_{0}\in\mathcal{P}(S)$,
$$
 Af(x)=\frac{\theta}{2}\int_{S}(f(y)-f(x))\nu_{0}(dy)
$$
describes the action of mutation and is called mutation operator.
 The domain of generator $\mathcal{L}$ is 
 $$
 \mathcal{D}(\mathcal{L})=\{\varphi\mid\varphi(\mu)=F(\langle f_{1},\mu\rangle,\cdots,\langle f_{k},\mu\rangle),F\in C^{2}(\mathbb{R}^{k}), f_{1},\cdots,f_{k}\in\mathcal{D}(A)\},
 $$
 where $\mathcal{D}(A)$ is the domain of mutation operator $A$.  Moreover, $\frac{\delta\varphi(\mu)}{\delta\mu(x)}$ can be regarded as directional derivative defined as 
$$
\frac{\delta\varphi(\mu)}{\delta\mu(x)}=\lim_{\epsilon\to0}\frac{\varphi(\mu+\epsilon\delta_{x})-\varphi(\mu)}{\epsilon}.
$$
  More background information can be found in \cite{MR2663265}.
 
 The Fleming-Viot process $Z_{t}$ is reversible(cf. \cite{MR1711591}) and has a unique stationary distribution (cf. \cite{MR1302692})
 $$
 \Pi_{\theta,\nu_{0}}(\cdot)=\mathbb{P}\left(\sum_{i=1}^{\infty}P_{i}(\theta)\delta_{\xi_{i}}\in\cdot\right),
 $$
where $(P_{1}(\theta),P_{2}(\theta),\cdots)$ follows the Poisson-Dirichlet distribution and $\{\xi_{i},i\geq1\}$, independent of $(P_{1}(\theta),P_{2}(\theta),\cdots)$, is an $i.i.d.$ sequence with common distribution $\nu_{0}$. The Poisson-Dirichlet distribution is first introduced by Kingman in \cite{MR0368264} and can be defined as the law of a descending order statistics of normalized jumps of gamma subordinator with L\'evy measure $\Lambda(dx)=\theta x^{-1}e^{-x}dx$. Also $\sum_{i=1}^{\infty}P_{i}(\theta)\delta_{\xi_{i}}$ is known as Dirichlet process, commonly used as a prior probability in nonparametric Bayesian statistics.
 
 Let $P(t,\mu,\cdot)$ be the transition probability of $Z_{t}$ starting with point $\mu$. Then $Z_{t}$ satisfies the following ergodic inequality(cf.\cite{MR1235429})
 \begin{equation}
 \|P(t,\mu,\cdot)-\Pi_{\theta,\nu_{0}}(\cdot)\|_{\mbox{var}}\leq(1+\theta)e^{-\frac{\theta}{2}t},\label{ergodic}
 \end{equation}
 where $\|\cdot\|_{\mbox{var}}$ is total variance norm.
 The parameter $\theta$ is usually defined as $\theta=4N_{e}\epsilon$, where $\epsilon$ is individual mutation rate and $N_{e}$ is effective population size. If $N_{e}\to\infty$, and $\epsilon$ is fixed, then $\theta\to\infty.$ Thus, for any fixed time $t$, however small,
 $$
 \lim_{\theta\to\infty}\|P(t,\mu,\cdot)-\Pi_{\theta,\nu_{0}}(\cdot)\|_{\mbox{var}}=0 
$$
 due to (\ref{ergodic}). It basically says that when the effective population size is large, the population will reach its equilibrium very quickly. Therefore, the small-time behaviour is quite critical for the comprehension of the whole evolution process. Furthermore, the time scale of evolution process can be quite large compared with lifespan of human beings. Hence, the observations obtained by biologists are of small-time scale. From this perspective, it is also interesting to look at the small-time behaviour.
 
 To understand the small-time behaviour of $Z_{t}$, we consider a time scale $t(\theta)$ approaching $0$ as $\theta\to\infty$. We are going to consider the asymptotic behaviour of $Z_{t(\theta)}$ when $\theta\to\infty$. When $t(\theta)=\frac{t}{\theta}$, the Wentzell-Freidlin large deviations for $Z_{t/\theta}$ have already been established in \cite{MR1815182},\cite{MR1649005} and \cite{MR1887170}. If $t(\theta)=a(\theta)t$ where $\lim_{\theta\to\infty}\theta a(\theta)=0$, the Schilder-type large deviations for $Z_{a(\theta)t}$ have been obtained in \cite{MR2184086}. Thus, under two typical small time scales $t(\theta)=\frac{t}{\theta}$ and $t(\theta)=a(\theta)t$, the asymptotic behaviour of $Z_{t(\theta)}$ is quite clear. Clearly, the small-time scale $a(\theta)t$ is smaller than $\frac{t}{\theta}$. The case where $t(\theta)$ is greater than $\frac{t}{\theta}$, however, is still untouched. Thus, our motivation here is to understand the small-time behaviour of $Z_{t(\theta)}$ with $\lim_{\theta\to\infty}\theta t(\theta)=\infty$ and $\lim_{\theta\to\infty}t(\theta)=0$. To this end, the weak limit of $Z_{t(\theta)}$ as $\theta\to\infty$ is obtained. One thing worthy of notice is that $\lim_{\theta\to\infty}t(\theta)=0$ and $t(\theta)$ is not necessarily confined with two aforementioned scales $t(\theta)=\frac{t}{\theta}$ and $t(\theta)=a(\theta)t$. To study the asymptotic behaviour of $Z_{t(\theta)}$ when $\lim_{\theta\to\infty}\theta t(\theta)=\infty$ and $\lim_{\theta\to\infty}t(\theta)=0$, we  consider large deviations for the transient sampling distributions of $Z_{t(\theta)}$. A simple reason for this treatment is that biologists only know the clustering property of a random sample from $Z_{t(\theta)}$. By a sample from $Z_{t(\theta)}$, we mean a sequence of $S$-valued random variables $\{X_{i},i\geq1\}$ such that, for any $n\geq1$, 
$$
\mathbb{P}(X_{1}\in dx_{1},\cdots, X_{n}\in dx_{n}|Z_{t(\theta)}=\mu)=\mu(dx_{1})\times\cdots\times \mu(dx_{n}).
$$
 Therefore, $\{X_{i},i\geq1\}$ is exchangeable and a reminiscent of paintbox process proposed by Kingman in \cite{MR509954}. We can therefore categorize a random sample $\{X_{1},\cdots,X_{n}\}$ by their allele types. Let $\{\mathcal{A}_{1}^n,\cdots,\mathcal{A}_{k}^n\}$ be the partitioned subclasses arranged decreasingly by class sizes. If we only care about species abundance, then we should turn to integer partition
$
(\#\mathcal{A}_{1}^n,\cdots,\#\mathcal{A}_{k}^n),
 $
 where $\#\mathcal{A}_{i}^n$ is the class size of $\mathcal{A}_{i}^n$. 
 
 By an integer partition, we mean an integer vector $\eta=(\eta_{1},\cdots,\eta_{l})$, where $\eta_{1}\geq\eta_{2}\geq\cdots\geq\eta_{l}$. Throughout this paper, $l(\eta)$ denotes the number of components of an integer partition $\eta$ and $|\eta|=\sum_{i=1}^{l(\eta)}\eta_{i}$. Define 
  $$
  \alpha_{i}(\eta)=\#\{j\mid \eta_{j}=i, 1\leq j\leq l(\eta)\}, 1\leq i\leq |\eta|.
  $$
  Then $(\alpha_{1}(\eta),\cdots,\alpha_{|\eta|}(\eta))$ is another notation for integer partition $\eta$. Obviously, $l(\eta)=\sum_{i=1}^{|\eta|}\alpha_{i}(\eta)$ and $|\eta|=\sum_{i=1}^{|\eta|}i\alpha_{i}(\eta).$  One may define an order $``<"$ in the partition set $\{\eta\mid |\eta|\geq1\}$ as follows:  $\eta<\xi$ if $|\eta|<|\xi|$ or $|\eta|=|\xi|$ but the first non-zero component of $\eta-\xi$ is negative.  In addition, two integer partition sets $\mathcal{J}^{o}=\{\eta\mid \eta_{l(\eta)}\geq2\}$ and $\mathcal{J}=\{(1), \eta\mid \eta_{l(\eta)}\geq2\}$ will be used repeatedly in this paper. One may project any partition $\eta=(\eta_{1},\cdots,\eta_{l})$ onto $\mathcal{J}$ as follows:
  $$
  \eta \mapsto\begin{cases}
  (1), &\mbox{ if } l=\alpha_{1}(\eta)\\
  (\eta_{1},\cdots,\eta_{l-\alpha_{1}(\eta)}) &\mbox{ if } l>\alpha_{1}(\eta).
  \end{cases}
  $$

 The distribution 
\begin{equation}
 P_{n}^{\theta}(\eta)=\mathbb{P}((\#\mathcal{A}_{1}^n,\cdots,\#\mathcal{A}_{k}^n)=\eta),\label{transient_sampling}
 \end{equation}
  is called sampling distribution. We therefore have a sequence of sampling distributions $\{P_{n}^{\theta},n\geq1\}$ derived from $\{X_{i},i\geq1\}$. For any positive integer $n$, since $\{X_{1},\cdots,X_{n}\}$ is a subsample of $\{X_{1},\cdots,X_{n},X_{n+1}\}$, we know $\{P_{n}^{\theta},n\geq1\}$ satisfies the consistency condition of Kolmogorov's extension theorem. A sequence of sampling distributions like $\{P_{n}^{\theta},n\geq1\}$ is called a partition structure introduced in \cite{MR0526801}. One important partition structure is the Ewens partition structure $\{E_{n}^{\theta},n\geq1\}$ determined by Dirichlet process $\sum_{i=1}^{\infty}P_{i}(\theta)\delta_{\xi_{i}}$. In this paper, when sample size $n$ is fixed, the large deviations for sampling distribution $\{P_{n}^{\theta},\theta>0\}$ is established. The large deviations for Ewens sampling distribution have already been obtained in \cite{MR2358634}. Compared with large deviations for Ewens sampling distribution, the large deviation results in this paper reveals a new structure which indicates a phase transition and is beyond expectation.
  
  This paper will proceed as follows. In section 2, the sampling distributions are discussed in detail. In section 3, the weak limit of $Z_{t(\theta)}$, where $\lim_{\theta\to\infty}\theta t(\theta)=+\infty,$ is obtained. In section 4, the large deviations for sampling distributions are established. The notation $``\sim"$ is used repeatedly in this paper. Here and after, we say $a(\theta)\sim b(\theta)$ if $\lim_{\theta\to\infty}\frac{a(\theta)}{b(\theta)}$.
  
 \section{Sampling Distribution}
  
  Since the probability measures on $S$ consist of two portions, discrete part and continuous part,  we can express any probability measure $\mu\in\mathcal{P}(S)$ as 
  $$
  \mu=\sum_{i=1}^{\infty}p_{i}\delta_{y_{i}}+\left(1-\sum_{i=1}^{\infty}p_{i}\right)\nu,
  $$
  where $\nu$ is a diffuse probability measure, $(p_{1},p_{2},\cdots)$ and $(y_{1},y_{2},\cdots)$ are masses and positions of the discrete part respectively. Define $D(\mu)$ to be $(p_{(1)},p_{(2)},\cdots)$, the descending order statistics of atomic masses of $\mu$. 
  
  Suppose that $X=\{X_{1},\cdots, X_{n}\}$ is a random sample from a population with distribution $\mu$ and let $D(\mu)$ be $x=(x_{1},x_{2},\cdots)$. Then the sampling distribution derived from the random sample $X$ can be calculated explicitly. To this end, we define the following mappings:
  \begin{itemize}
  \item[1,] $G: X\to (\mathcal{A}_{1}^n,\cdots,\mathcal{A}_{k}^n)$ which maps a random sample to its partition.
  \item[2,] $T:X\to(\#\mathcal{A}_{1}^n,\cdots,\#\mathcal{A}_{k}^n)$ which maps a random sample to an integer partition.
  \end{itemize}
 In particular, $T(X)=(n)$ means that all sampled individuals are of the same type. So
  $$
   \mathbb{P}( T(X)=(n))=P(X_{1}=\cdots=X_{n})=\varphi_{n}(D(\mu))=\sum_{i=1}^{\infty}x_{i}^{n}.
  $$
In general, for $\eta=(\eta_{1},\cdots,\eta_{l})$, $\{T(X)=\eta\}$ means all the random samples that can be classified as the integer partition $\eta$. One specific random sample satisfies this condition is the partition $\mathcal{A}^{*}(n)=(\mathcal{A}_{1}^n,\cdots,\mathcal{A}_{l}^n)$, where 
 \begin{equation}
 \mathcal{A}_{i}^n=\left\{X_{\sum_{j=1}^{i-1}\eta_{j}+1},\cdots,X_{\sum_{j=1}^{i}\eta_{j}}\right\},1\leq i\leq l,\label{P1}
 \end{equation}
  where  all individuals in $\mathcal{A}_{i}^n$ share the same type.
 Due to the exchangeability of sample $X$, it is easy to show that 
$$
\mathbb{P}(T(X)=\eta)=\frac{n!}{\prod_{i=1}^{l}\eta_{i}!\prod_{i=1}^{n}(i!)^{\alpha_{i}}}p_{\eta}^{o}(x),  
  $$
  where $p_{\eta}^{o}(x)=\mathbb{P}(G(X)=\mathcal{A}^{*}(n))$.  Consider a subsample 
  $$
  X_{sub}=\left\{X_{1},\cdots,X_{\sum_{i=1}^{l-1}\eta_{j}}\right\}
  $$ 
  of the random sample $X$. Similarly, one particular partition of  $X_{sub}$ is $\mathcal{A}^{*}(n-\eta_{l})=(\mathcal{A}_{1}^{n-\eta_{l}},\cdots,\mathcal{A}_{l-1}^{n-\eta_{l}})$, where
  \begin{equation}
 \mathcal{A}_{i}^{n-\eta_{l}}=\left\{X_{\sum_{j=1}^{i-1}\eta_{j}+1},\cdots,X_{\sum_{j=1}^{i}\eta_{j}}\right\},1\leq i\leq l-1,\label{P2}
\end{equation} 

Then $ p_{(\eta_{1},\cdots,\eta_{l-1})}^{o}(x)=\mathbb{P}(G(X_{sub})=\mathcal{A}^{*}(n-\eta_{l}))$. We claim that  
\begin{equation}
p_{\eta}^{o}(x)=p^{o}_{(\eta_{1},\cdots,\eta_{l-1})}(x)\varphi_{\eta_{l}}(x)-\sum_{j=1}^{l-1}p^{o}_{\eta^{j}}(x),\label{iteration}
\end{equation}
Indeed,
  \begin{align*}
  &p^{o}_{(\eta_{1},\cdots,\eta_{l-1})}(x)\varphi_{\eta_{l}}(x)\\
  =&
  \mathbb{P}(\{G(X_{sub})=\mathcal{A}^{*}(n-\eta_{l})\}\cap \{X_{\sum_{i=1}^{l-1}\eta_{i}+1}=\cdots=X_{n}\} )\\
=&\sum_{i=1}^{l-1}\mathbb{P}(\{G(X_{sub})=\mathcal{A}^{*}(n-\eta_{l})\}\cap\{X_{\sum_{i=1}^{l-1}\eta_{i}+1}=\cdots=X_{n} \mbox{ share type in }\mathcal{A}_{i}^{n-\eta_{l}}\})\\
&+ \mathbb{P}(\{G(X_{sub})=\mathcal{A}^{*}(n-\eta_{l})\}\cap \{X_{\sum_{i=1}^{l-1}\eta_{i}+1}=\cdots=X_{n} \mbox{ have brand new type }\})\\
=&\sum_{j=1}^{l-1}p^{o}_{\eta^{j}}(x)+p_{\eta}^{o}(x)
 \end{align*}
where $\eta^{j}$ is a decreasing arrangement of $(\eta_{1}^{j},\cdots,\eta_{l-1}^{j})$ and 
$$
\eta^{j}_{i}=\begin{cases}
\eta_{i}+\eta_{l},& \mbox{ if } i=j\\
\eta_{i},& \mbox{ otherwise }. 
\end{cases}
$$
 In particular, when $\eta=(1)$, $p^{o}_{\eta}(x)=\varphi_{1}(x)=1$.  Making use of (\ref{iteration}), we can express $p_{\eta}^{o}(x)$ in terms of $\varphi_{k}(x),k\geq1.$
 
  \begin{pro}\label{partition_representation}
  Denote $\pi(l,d), 1\leq d\leq l,$ to be the set of partitions of $$\{1,2,\cdots,l\}$$ into $d$ subclasses $(\beta_{1},\cdots,\beta_{d})$, satisfying
  $$
 \min\beta_{1}<\min\beta_{2}<\cdots<\min\beta_{d}.
 $$
 Let $|\beta_{i}|$ be the cardinality of $\beta_{i}$. Then,  for a given integer partition $\eta=(\eta_{1},\cdots,\eta_{l})$,  \begin{equation}
 p_{\eta}^{o}(x)=\sum_{d=1}^{l}(-1)^{l-d}\sum_{\beta\in\pi(l,d)}(|\beta_{1}|-1)!\cdots (|\beta_{d}|-1)!
 \varphi_{\sum_{i\in\beta_{1}}\eta_{i}}(x)\cdots\varphi_{\sum_{i\in\beta_{d}}\eta_{i}}(x).\label{new_function}
  \end{equation}
  \end{pro}
  \begin{rmk}
  \item[1,] 
  If $\mu$ is a purely atomic probability measure, then $\sum_{i=1}^{\infty}x_{i}=1$ and we claim $p_{\eta}^{o}(x)$  is equal to
  \begin{equation}
  q^{o}_{\eta}(x)=\sum_{i_{1},\cdots,i_{l}\neq} x_{i_{1}}^{\eta_{1}}\cdots x_{i_{l}}^{\eta_{l}}\label{extension}
  \end{equation}
  Indeed,  it is not difficult to show that $q_{\eta}^{o}(x)$ satisfies condition $(\ref{iteration})$, $q_{(n)}^{o}(x)=\varphi_{n}(x)$ and $q_{(1)}^{o}(x)=1$. Running the similar argument in Proposition \ref{partition_representation} yields $q_{\eta}^{o}(x)=p_{\eta}^{o}(x).$
\item[2,]  Consider 
  $$
  \triangledown_{\infty}=\left\{x\in[0,1]^{\infty}\Big| x_{1}\geq x_{2}\cdots\geq0,\sum_{i=1}^{\infty}x_{i}=1 \right\}
  $$
  and
  $$
  \bar{\triangledown}_{\infty}=\left\{x\in[0,1]^{\infty}\Big| x_{1}\geq x_{2}\cdots\geq0,\sum_{i=1}^{\infty}x_{i}\leq1 \right\}. 
  $$
  Then $\triangledown_{\infty}$ is a dense subspace of $\bar{\triangledown}_{\infty}$. Under product topology, $\bar{\triangledown}_{\infty}$ is compact. It is easy to see that $p_{\eta}^{o}$ is a continuous extension of $(\ref{extension})$.
 \end{rmk}
  Now we can easily obtain the following sampling distributions.
  \begin{theorem}\label{sampling_distribution_theorem}
  Sampling distributions $\{P_{n}^{\theta},n\geq1\}$ of a random sample from $Z_{t(\theta)}$ is
  $$
  P_{n}^{\theta}(\eta)=\frac{n!}{\prod_{i=1}^{l}\eta_{i}!\prod_{i=1}^{n}(i!)^{\alpha_{i}(\eta)}}\mathbb{E}p_{\eta}^{o}(D(Z_{t(\theta)}))
  $$
  \end{theorem}
\begin{rmk}
\item[1,] The sampling distribution $\{P_{n}^{\theta}(\eta),n\geq1\}$ can also be represented as
 \begin{equation}
 P_{n}^{\theta}(\eta)= \int_{\bar{\triangledown}_{\infty}} \frac{n!}{\prod_{i=1}^{l}\eta_{i}!\prod_{i=1}^{n}(i!)^{\alpha_{i}(\eta)}}p_{\eta}^{o}(x) \nu_{\theta}(dx),  \label{new_represen}
\end{equation}
 where $\nu_{\theta}(dx)$ is the distribution of a descending order statistics of atomic masses of $Z_{t(\theta)}$. Because $p_{\eta}^{o}(x)$ is continuous in $\bar{\triangledown}_{\infty}$, the weak convergence of $\nu_{\theta}$ leads to the point wise convergence of sampling distributions.
 \item[2,] The weak limit of $\nu_{\theta}$ is not necessarily a purely atomic probability measure. When $\mu$ is a diffuse probability measure, $P_{n}^{\theta}(\eta)=\delta_{(1,\cdots,1)}(\eta)$,i.e. all individuals have different types.
  \end{rmk}
  
  \begin{exm}
  Let $\{E^{\theta}_{n},n\geq1\}$ be a family of sampling distributions from Dirichlet process $\sum_{i=1}^{\infty}P_{i}(\theta)\delta_{\xi_{i}}$. Then, for a given partition $\eta=(\eta_{1},\cdots,\eta_{l})$,
  \begin{equation}
  E^{\theta}_{n}(\eta)=\frac{n!}{\prod_{i=1}^{l}\eta_{i}\prod_{i=1}^{n}(i!)^{\alpha_{i}(\eta)}}\frac{\theta^{l}}{\theta_{(n)}}, \label{Ewens_Sampling_Formula}
  \end{equation}
  which is the Ewens sampling formula(cf. \cite{MR2663265}).
 \end{exm}

 \section{Weak Limit of $Z_{t(\theta)}$}
 Let $t(\theta)$ be a small-time scale, i.e. $\lim_{\theta\to\infty}t(\theta)=0.$ We are going to pinpoint the weak limit of $Z_{t(\theta)}$. To this end, we denote $P(t,\mu,d\nu)$ to be the transition probability of $Z_{t}$. The explicit representation of $P(t,\mu,d\nu)$ obtained by Ethier and Griffiths in \cite{MR1235429} is the following:
 \begin{align*}
 P(t,\mu,d\nu)=&d_{0}^{\theta}(t)\Pi_{\theta,\nu_{0}}(d\nu)\\
 &+\sum_{n=1}^{\infty}d_{n}^{\theta}(t)\int_{S^n}\mu^n(dx_{1}\times\cdots\times dx_{n})\Pi_{n+\theta,(n+\theta)^{-1}(\sum_{i=1}^{n}\delta_{x_{i}}+\theta\nu_{0})}(d\nu),
 \end{align*}
 where
 \begin{align*}
 d_{0}^{\theta}(t)&=1-\sum_{m=1}^{\infty}\frac{2m-1+\theta}{m!}(-1)^{m-1}\theta_{(m-1)}e^{-\lambda_{m}t}\\
 d_{n}^{\theta}(t)&=\sum_{m=n}^{\infty}\frac{2m-1+\theta}{m!}(-1)^{m-n}\binom{m}{n}(n+\theta)_{(m-1)}e^{-\lambda_{m}t}, n\geq1.
 \end{align*}
 $\lambda_{0}=0, \lambda_{m}=\frac{m(m-1+\theta)}{2},m\geq1$. $\{0,-\lambda_{m},m\geq1\}$ are the eigenvalues of $Z_{t}.$ Moreover, $\Pi_{n+\theta,(n+\theta)^{-1}(\sum_{i=1}^{n}\delta_{x_{i}}+\theta\nu_{0})}$ is the distribution of a general Dirichlet process 
 $
 \sum_{i=1}^{\infty}P_{i}(n+\theta)\delta_{\xi_{i}}
 $ where $\xi_{i}$'s are $i.i.d.$ $\frac{n}{n+\theta}\sum_{i=1}^{n}\delta_{x_{i}}+\frac{\theta}{n+\theta}\nu_{0}$, and $(P_{1}(n+\theta),\cdots)$ follows Poisson-Dirichlet distribution $\mbox{PD}(n+\theta)$.

  The following lemma is from \cite{MR1235429} as well.
 \begin{lemma}\label{moments_formula}
 Let $f_{1},f_{2},\cdots,f_{m}\in C(S)$. Define  $a_{[n]}=a(a-1)\cdots (a-n+1)$, and $a_{(n)}=a(a+1)\cdots (a+n-1)$. Then an explicit formula for 
$$
\int_{\mathcal{P}(S)} \langle f_{1},\nu\rangle\cdots  \langle f_{m},\nu\rangle P(t,\mu,d\nu),
$$ 
is 
 \begin{align*}
\sum_{n=0}^{\infty}d_{n}^{\theta}(t)\sum_{M\subset\{1,\cdots,m\}}&\frac{1}{(n+\theta)_{(m)}}
\left\{\sum_{k=1}^{|M|}n_{[k]}\sum_{\beta\in\pi(|M|,k)}|\beta_{1}|!\cdots|\beta_{k}|!\prod_{j=1}^{k}\langle\prod_{i\in\beta_{j}}f_{i},\mu \rangle\right\}\\
&\times\left\{\sum_{l=1}^{|M^{c}|}\sum_{\gamma\in\pi(|M^{c}|,l)}(|\gamma_{1}|-1)!\cdots(|\gamma_{l}|-1)!\theta^{l}\prod_{j=1}^{l}\langle \prod_{i\in\gamma_{j}}f_{i},\nu_{0}\rangle\right\}.
 \end{align*}
 \end{lemma}

 Define $\alpha_{m}^{k}(t), 0\leq k\leq m,$ to be the coefficient of the term
 $$
 \left\{|\beta_{1}|!\cdots|\beta_{k}|!\prod_{j=1}^{k}\langle\prod_{i\in\beta_{j}}f_{i},\mu \rangle\right\}\times \left\{(|\gamma_{1}|-1)!\cdots(|\gamma_{l}|-1)!\theta^{l}\prod_{j=1}^{l}\langle \prod_{i\in\gamma_{j}}f_{i},\nu_{0}\rangle \right\},
 $$
 i.e.
 $$
 \alpha_{m}^{k}(t)=\sum_{n=0}^{\infty}\frac{n_{[k]}}{(n+\theta)_{(m)}}d_{n}^{\theta}(t).
 $$
 
 By (3.15) in \cite{MR1235429}, we know that $\alpha_{m}^{k}(t)$ satisfies the following differential equations:
 \begin{align*}
 \frac{d \alpha_{m}^{k}(t)}{dt}=-\lambda_{m}\alpha_{m}^{k}(t)+\frac{m-k}{2}\alpha_{m-1}^{k}(t)
 \end{align*}
 with initial conditions $\lim_{t\to0}\alpha_{m}^{k}(t)=\delta_{k,m}(0\leq k\leq m)$.
 Therefore,
 \begin{align*}
\alpha_{m}^{k}(t)=\frac{(m-k)!}{2^{m-k}}e^{-\lambda_{m}t}&\idotsint\limits_{0\leq s_{k}\leq\cdots\leq s_{m-1}\leq t}\\
&e^{(\lambda_{m}-\lambda_{m-1})s_{m-1}}\cdots e^{(\lambda_{k+1}-\lambda_{k})s_{k}}ds_{m-1}\cdots ds_{k}.
 \end{align*}
We consider a small-time scale $t(\theta)$, such that $\lim_{\theta\to\infty}\theta t(\theta)=c\in(0,\infty)$. Under transform 
 $$
 s_{i}=\frac{u_{i}}{\theta}, k\leq i\leq m-1,
 $$
 \begin{align*}
 \alpha_{m}^{k}(t(\theta))=&\frac{(m-k)!}{\theta^{m-k}2^{m-k}}e^{-\frac{\lambda_{m}}{\theta}\theta t(\theta)}\idotsint\limits_{0\leq u_{k}\leq\cdots\leq u_{m-1}\leq \theta t(\theta)}\\
 &e^{\frac{(\lambda_{m}-\lambda_{m-1})}{\theta}u_{m-1}}\cdots e^{\frac{(\lambda_{k+1}-\lambda_{k})}{\theta}u_{k}}du_{m-1}\cdots du_{k}\\
\sim&\frac{(m-k)!}{\theta^{m-k}2^{m-k}}e^{-\frac{mc}{2}}\idotsint\limits_{0\leq u_{k}\leq\cdots\leq u_{m-1}\leq c}e^{\frac{u_{m-1}}{2}}\cdots e^{\frac{u_{k}}{2}}du_{m-1}\cdots du_{k} \\
=& \frac{1}{\theta^{m-k}2^{m-k}}e^{-\frac{mc}{2}}\left(\int_{0}^{c}e^{\frac{s}{2}}ds\right)^{m-k}\\
=&\frac{1}{\theta^{m-k}}e^{-\frac{kc}{2}}(1-e^{-\frac{c}{2}})^{m-k}.
 \end{align*}
Thus,
\begin{equation}
\alpha_{m}^{k}(t(\theta))\sim\frac{1}{\theta^{m-k}}e^{-\frac{kc}{2}}(1-e^{-\frac{c}{2}})^{m-k}.\label{limit}
\end{equation} 
If $\lim_{\theta\to\infty}\theta t(\theta)=0$ or $\lim_{\theta\to\infty}\theta t(\theta)=0$, however, we have 
\begin{equation}
\lim_{\theta\to\infty}\theta^{m-k}\alpha_{m}^{k}(t(\theta))=0, \label{lower_order}
\end{equation}
and 
\begin{equation}
\lim_{\theta\to\infty}\theta^{m-k}\alpha_{m}^{k}(t(\theta))=\infty.
 \label{lower_order}
\end{equation}
 \begin{theorem} \label{weak}
 When $\lim_{\theta\to\infty}t(\theta)=0$ and $\lim_{\theta\to\infty}\theta t(\theta)=c$ where $c$ can be $0$, $\infty$ and any positive number,  then $Z_{t(\theta)}$ converges to 
 $$
 m_{c}=e^{-\frac{c}{2}}\mu+(1-e^{-\frac{c}{2}})\nu_{0}
 $$
 in probability.
 \end{theorem}
 \begin{rmk}
\item[1,] When $t(\theta)=\frac{t}{\theta}$ and $t(\theta)=a(\theta)t$ where $\lim_{\theta\to\infty}\theta a(\theta)=0$, this result can also be derived from path-level large deviations in \cite{MR1815182},\cite{MR1649005}, \cite{MR1887170},
\cite{MR2184086}.
\item[2,] When $\lim_{\theta\to\infty}t(\theta)=0$ and $\lim_{\theta\to\infty}\theta t(\theta)=\infty$, this result is new. Also the proof of Theorem \ref{weak} is necessary because when $\lim_{\theta\to\infty}(\theta+1)e^{-\theta t(\theta)}\neq0$ we can not use the ergodic inequality of $Z_{t}$ to find the weak limit of $Z_{t(\theta)}$.
\item[3,] Because $m_{c}$ is a single point in $\mathcal{P}(S)$, we can show that $\mathbb{P}(Z_{t(\theta)}\in\cdot)$ converges to $\delta_{m_{c}}$ weakly.
\end{rmk}
 \begin{proof}
  For any $f\in C(S)$, we claim that the variance of $Z_{t(\theta)}$vanishes as $\theta\to\infty$ and 
  $$
  \lim_{\theta\to\infty}\mathbb{E}\langle Z_{t(\theta)},f_{i}\rangle=\langle m_{c},f_{i}\rangle.
  $$ 
 Thus for any $\epsilon>0, k>1$ satisfying $2^{1-k}<\epsilon$,
 \begin{align*}
 \mathbb{P}(d_{w}(Z_{t(\theta)},m_{c})>\epsilon)\leq&\sum_{i=1}^{k}\mathbb{P}\left(|\langle Z_{t(\theta)}-m_{c},f_{i}\rangle|\geq\frac{2^{i-1}\epsilon}{k}\right)\\
 \leq&\sum_{i=1}^{k}\left(\frac{k}{2^{i-1}\epsilon}\right)^{2}\left[\mbox{Var}(\langle f_{i},Z_{t(\theta)}\rangle)+(\mathbb{E}\langle Z_{t(\theta)},f_{i}\rangle-\langle m_{c},f_{i}\rangle)^{2}\right]\\
 \to&0.
 \end{align*}
 The proof is thus completed.  Now we are going to show the previous claims.

  By Lemma  \ref{moments_formula} and (\ref{limit}), whatever $c$ may be, we have
 \begin{align*}
 \mathbb{E}\langle f,Z_{t(\theta)}\rangle=&e^{-\frac{\theta}{2}t(\theta)}\langle f,\mu\rangle+\left(1-e^{-\frac{\theta}{2}t(\theta)}\right)\langle f,\nu_{0}\rangle\\
 \to& e^{-\frac{c}{2}}\langle f,\mu\rangle+(1-e^{-\frac{c}{2}})\langle f,\nu_{0}\rangle.\\
 \mathbb{E}\langle f,Z_{t(\theta)}\rangle^{2}=&\alpha_{2}^{1}(t(\theta))2!\langle f^{2},\mu\rangle+\alpha_{2}^{2}(t(\theta))\langle f,\mu\rangle^{2}+2\alpha_{2}^{1}(t(\theta))\theta\langle f,\mu\rangle\langle f,\nu_{0}\rangle\\
 &+\alpha_{2}^{0}(t(\theta))\theta^{2}\langle f,\nu_{0}\rangle^{2} +\alpha_{2}^{0}\theta\langle f^{2},\nu_{0}\rangle\\
 \to& e^{-c}\langle f,\mu\rangle^{2}+2e^{-\frac{c}{2}}(1-e^{-\frac{c}{2}})\langle f,\mu\rangle\langle f,\nu_{0}\rangle+(1-e^{-\frac{c}{2}})^{2}\langle f,\nu_{0}\rangle^{2}\\
 =&(e^{-\frac{c}{2}}\langle f,\mu\rangle+(1-e^{-\frac{c}{2}})\langle f,\nu_{0}\rangle)^{2}.
\end{align*}
Therefore, as $\theta\to\infty,$
$$
\mbox{Var}(\langle f,Z_{t(\theta)}\rangle)\to0.
$$
 \end{proof}
 
 Generally, the weak convergence of the marginal laws of $Z_{t(\theta)}$ might not necessarily guarantee the weak convergence of the distributions of $D(Z_{t(\theta)})$. Therefore, based on the weak convergence of the marginal laws of $Z_{t(\theta)}$, we can not deduce the point wise convergence of the sampling distributions $\{P_{n}^{\theta},\theta>0\}$ derived from a random sample from $Z_{t(\theta)}$. Ethier and Kurtz in \cite{MR1302692} proposed  weak atomic topology, generated by metric
 \begin{align*}
 d_{a}(\mu,\nu)=d_{w}(\mu,\nu)&+\sup_{0\leq\epsilon\leq1}\Bigg|\int_{S}\int_{S}\left(1-\frac{r(x,y)}{\epsilon}\right)_{+}\mu(dx)\mu(dy)\\
 &-\int_{S}\int_{S}\left(1-\frac{r(x,y)}{\epsilon}\right)_{+}\nu(dx)\nu(dy)\Bigg|,
 \end{align*}  
 where $r(x,y)$ is the metric function in type space $S$ and 
 $$
 \left(1-\frac{r(x,y)}{\epsilon}\right)_{+}=\begin{cases}
 1-\frac{r(x,y)}{\epsilon}, & \mbox{ if } r(x,y)\leq\epsilon\\
 0, & \mbox{ ohterwise }.
 \end{cases}
 $$ Obviously, the weak atomic topology is stronger than weak topology, it can guarantees the weak convergence of $D(Z_{t(\theta)})$ if the law of $Z_{t(\theta)}$ converges in weak atomic topology.
 \begin{cor}\label{cor}
 Under conditions of Theorem \ref{weak}, $Z_{t(\theta)}$ also converges to
 $m_{c}$ in probability in $(\mathcal{P}(S),d_{a}).$ 
 \end{cor}
 
 In fact, define $X_{t}=D(Z_{t})$. Then $X_{t}$ is an atomic diffusion constructed in \cite{MR615945}. The generator of $X_{t}$ is 
 $$
 G=\frac{1}{2}\sum_{i=1}^{\infty}x_{i}(\delta_{ij}-x_{j})\frac{\partial^{2}}{\partial x_{i}\partial x_{j}}-\frac{\theta}{2}\sum_{i=1}^{\infty}x_{i}\frac{\partial}{\partial x_{i}}.
 $$ 
 In \cite{MR615945}, the Poisson-Dirichlet distribution $\mbox{PD}(\theta)$ was proved to be the stationary distribution of $X_{t}$. Either also obtained an explicit transition density $p(t,x,y)$ of $X_{t}$ in \cite{MR1174426}, where
 \begin{equation}
 p(t,x,y)=1+\sum_{m=2}^{\infty}e^{-\lambda_{m}t}Q_{m}(x,y).\label{density}
 \end{equation}
 Here $Q_{m}(x,y)=\sum_{|\eta|=m}\chi_{\eta}(x)\chi_{\eta}(y),$ where $\eta\in\mathcal{J}$ and $\chi_{\eta}(x)$ are orthonormal eigenfunctions associated with eigenvalues $-\lambda_{m}$ of generator $G$ in Hilbert space $L^2(\mbox{PD}(\theta))$. For a given integer partition $\eta=(\eta_{1},\cdots,\eta_{l})$ and $\eta\in\mathcal{J}$, we define $\varphi_{\eta}(x)=\varphi_{\eta_{1}}(x)\cdots\varphi_{\eta_{l}}(x)$. Ethier concluded in  \cite{MR1174426} that all orthonormal eigenfunctions $\chi_{\eta}(x)$ can be obtained by applying Gram-Schmidt orthogonalization to 
 $
 \{\varphi_{\eta}(x) \mid\eta\in\mathcal{J}
\}.
 $ 
 
 Let $\langle f,g\rangle_{\theta}=\int_{\bar{\triangledown}_{\infty}}fd\mbox{PD}(\theta)$ be the inner product in $L^2(\mbox{PD}(\theta))$, and $\|f\|=\sqrt{\langle f,g\rangle_{\theta}}$ be the induced norm. If we adopt the order $``<"$, the  Gram-Schmidt orthogonalization process can go as follow
 $$
\psi^{\theta}_{1}(x)=1, \quad \psi_{2}^{\theta}(x)=\varphi_{2}(x)-\frac{1}{1+\theta}, 
$$
and generally
\begin{align}
\psi_{\eta}^{\theta}(x)=&\varphi_{\eta}-\sum_{\xi<\eta,\xi\in\mathcal{J}}\langle\varphi_{\eta},\psi_{\xi}^{\theta}\rangle_{\theta}\psi_{\xi}^{\theta}(x)\quad \eta\in\mathcal{J}.\label{gram_schmidt}
\end{align}  
 Define $\chi_{\eta}^{\theta}(x)=\frac{\psi_{\eta}^{\theta}(x)}{||\psi_{\eta}^{\theta}||},\eta\in\mathcal{J}.$ Then $\{\chi_{\eta}^{\theta}\mid\eta\in\mathcal{J}\}$ is the set of orthonormal eigenfunctions.  One can easy show the following lemma.
 \begin{lemma}\label{orthonormal}
 \begin{align}
\int_{\bar{\triangledown}_{\infty}}Q_{m}(x,y)\psi_{\eta}^{\theta}(y)\mbox{PD}(\theta)(dy)=\psi_{\eta}^{\theta}(x)\delta_{|\eta|,m}.\label{orthogonal_integration}
\end{align} 
\end{lemma}
To prove Corollary (\ref{cor}), we need the following lemma. 
\begin{lemma}\label{var}
As $\theta\to\infty$,
\begin{equation}
 \sum_{x\in S}Z_{t(\theta)}(\{x\})^{2}\to \sum_{x\in S}m_{c}(\{x\})^{2} \mbox{ in probability }.\label{atomic}\end{equation}  
 \end{lemma}
 \begin{proof}
 It suffices to show that $\mathbb{E} \sum_{x\in S}Z_{t(\theta)}(\{x\})^{2}\to\sum_{x\in S}m_{c}(\{x\})^{2}$ and the variance of $ \sum_{x\in S}Z_{t(\theta)}(\{x\})^{2}$ vanishes as $\theta\to\infty$. Since $ \sum_{x\in S}Z_{t(\theta)}(\{x\})^{2}=\varphi_{2}(D(Z_{t(\theta)}))=\varphi_{2}(X_{t(\theta)})$, we have
 $$
 Var\left(\sum_{x\in S}Z_{t(\theta)}(\{x\})^{2}\right)=\mathbb{E}\varphi_{2}^2(X_{t(\theta)})-[\mathbb{E}\varphi_{2}(X_{t(\theta)})]^2.
 $$
 Making use of the transition density of $X_{t}$ in (\ref{density}), we have
 \begin{align*}
 \mathbb{E}\varphi_{2}(X_{t(\theta)})=&\int_{\bar{\triangledown}_{\infty}}\varphi_{2}(y)p(t(\theta),x,y)\mbox{PD}(\theta)(dy)\\
 =& \int_{\bar{\triangledown}_{\infty}}\varphi_{2}(y)\mbox{PD}(\theta)(dy)+\sum_{m=2}^{\infty}e^{-\lambda_{m}t(\theta)}\int_{\bar{\triangledown}_{\infty}}\varphi_{2}(y)Q_{m}(x,y)\mbox{PD}(\theta)(dy)\\
  \mathbb{E}\varphi_{2}^2(X_{t(\theta)})=&\int_{\bar{\triangledown}_{\infty}}\varphi_{2}^2(y)p(t(\theta),x,y)\mbox{PD}(\theta)(dy)\\
 =& \int_{\bar{\triangledown}_{\infty}}\varphi_{2}^2(y)\mbox{PD}(\theta)(dy)+\sum_{m=2}^{\infty}e^{-\lambda_{m}t(\theta)}\int_{\bar{\triangledown}_{\infty}}\varphi_{2}^2(y)Q_{m}(x,y)\mbox{PD}(\theta)(dy)  
 \end{align*}
 Due to Gram-Schmidt orthogonalization process, we have 
 \begin{align}
 \varphi_{2}(x)=&\psi_{2}^{\theta}(x)+\frac{1}{1+\theta} \label{1}\\
 \varphi_{2}^2(x)=&\psi_{2,2}^{\theta}(x)+\langle\varphi_{2}^2,\psi_{3}^{\theta}\rangle_{\theta}\psi_{3}^\theta(x)+\langle\varphi_{2}^2,\psi_{2}^{\theta}\rangle_{\theta}\psi_{2}^\theta(x)+\langle\varphi_{2}^2,1\rangle_{\theta}. \label{2}
 \end{align}
 Replacing $\varphi_{2},\varphi_{2}^2$ in $\mathbb{E}\varphi_{2}(X_{t(\theta)})$ and $\mathbb{E}\varphi_{2}^2(X_{t(\theta)})$ by  (\ref{1}) and (\ref{2}), we have
 $$
 \mathbb{E}\varphi_{2}(X_{t(\theta)})=\frac{1}{1+\theta}+e^{-(1+\theta)t(\theta)}\psi_{2}^{\theta}(x)=\frac{1}{1+\theta}+e^{-(1+\theta)t(\theta)}\left(\varphi_{2}(x)-\frac{1}{1+\theta}\right) $$
 and
 \begin{align*}
  \mathbb{E}\varphi_{2}^2(X_{t(\theta)})=&\langle\varphi_{2}^2,1\rangle_{\theta}+e^{-\lambda_{4}t(\theta)}\psi_{2,2}^{\theta}(x)+e^{-\lambda_{3}t(\theta)}\langle\varphi_{2}^2,\psi_{3}^{\theta}\rangle_{\theta}\psi_{3}^\theta(x)\\
  &+e^{-\lambda_{2}t(\theta)}\langle\varphi_{2}^2,\psi_{2}^{\theta}\rangle_{\theta}\psi_{2}^\theta(x)  
  \end{align*}
  due to Lemma \ref{orthonormal}. Since $\int_{\bar{\triangledown}_{\infty}}\varphi_{2}d\mbox{PD}(\theta)=\frac{1}{1+\theta}\to0$ as $\theta\to\infty$, then we know $\mbox{PD}(\theta)$ converges weakly to $(0,0,\cdots)$ as $\theta\to\infty$. Therefore, 
  $$
 \lim_{\theta\to\infty} \mathbb{E} \sum_{x\in S}Z_{t(\theta)}(\{x\})^{2}=\lim_{\theta\to\infty} \mathbb{E}\varphi_{2}(X_{t(\theta)})=e^{-c}\varphi_{2}(x)=\sum_{x\in S}m_{c}(\{x\})^{2}  
 $$
 and 
 $$
 \lim_{\theta\to\infty} \mathbb{E}\varphi_{2}^2(X_{t(\theta)})=e^{-2c}\varphi_{2}^2(x)=\left[\sum_{x\in S}m_{c}(\{x\})^{2}\right]^2.
 $$
  The proof is thus completed.
\end{proof}
 [\textbf{Proof of Corollary \ref{cor}}]:
 
  \begin{proof}
 First of all, $Z_{t(\theta)}\to m_{c}$ in probability in $(\mathcal{P},d_{a})$ is equivalent to that, for any sequence $\{\theta^{'},\theta^{'}>0\}$, there is a subsequence $\{\theta_{n}^{'}\}$ such that $Z_{t(\theta_{n}^{'})}\to m_{c}$ almost surely in $(\mathcal{P},d_{a})$. Therefore, for a given sequence $\{\theta^{'}, \theta^{'}>0\}$, it suffices to find a subsequence $\{\theta_{n},n\geq1\}$ such that $Z_{t(\theta_{n})}\to m_{c}$ almost surely in $(\mathcal{P},d_{a})$. To this end,
  since $Z_{t(\theta^{'})}\to m_{c}$ in probability in $(\mathcal{P},d_{w})$, there exists a subsequence $\{\theta_{n}^{'},n\geq1\}$ such that $Z_{t(\theta_{n}^{'})}\to m_{c}$ almost surely in $(\mathcal{P},d_{w})$. 
  
  By Lemma 4.2 in \cite{MR1205982}, we only need to find a subsequence $\{\theta_{n},n\geq1\}$ such that 
\begin{equation}
 \sum_{x\in S}Z_{t(\theta_{n})}(\{x\})^{2}\to \sum_{x\in S}m_{c}(\{x\})^{2} \mbox{ almost surely },\label{atomic}\end{equation}    
 By Lemma \ref{var}, we know 
 $$
 \sum_{x\in S}Z_{t(\theta_{n}^{'})}(\{x\})^{2}\to \sum_{x\in S}m_{c}(\{x\})^{2} \mbox{ in probability}.
  $$
 Therefore, we can find a subsequence $\{\theta_{n},n\geq1\}$ such that (\ref{atomic}) is true. Because of Lemma 4.2 in \cite{MR1205982}, $Z_{t(\theta_{n})}\to m_{c}$ almost surely in $(\mathcal{P}, d_{a})$.  Hence,
$Z_{t(\theta)}\to m_{c}$ in probability in $(\mathcal{P},d_{a})$. 
\end{proof}

 Because for any sequence $Z_{t(\theta_{n})}$ converging to $m_{c}$ in probability in $(\mathcal{P},d_{a})$, there is a subsequence $Z_{t(\theta_{n}^{'})}$ which converges almost surely in $(\mathcal{P},d_{a})$. So $D(Z_{t(\theta_{n}^{'})})$ converges almost surely to $D(m_{c})$. This indicates that $D(Z_{t(\theta_{n})})$ converges to $D(m_{c})$ in probability in $\bar{\triangledown}_{\infty}$. We can also show that $D(Z_{t(\theta)})$ converges weakly to $D(m_{c})$. Therefore, due to Theorem \ref{sampling_distribution_theorem}, the point wise convergence of sampling distributions is also true.

\begin{theorem}\label{asymptotic_sampling}
Let $P^{\theta}_{n}$ be the sampling distributions of $Z_{t(\theta)}$.
When $\lim_{\theta\to\infty}\theta t(\theta)=c,$ we have, for any integer partition $\eta$,  
$$
\lim_{\theta\to\infty}P_{n}^{\theta}(\eta)=\frac{n!}{\prod_{i=1}^{l}\eta_{i}!\prod_{i=1}^{n}(i!)^{\alpha_{i}}}p_{\eta}^{o}(D(m_{c})).
$$
In particular, when $c=\infty$,
$$
\lim_{\theta\to\infty}P_{n}^{\theta}(\eta)=\begin{cases}
1,& \mbox{ if } \eta=(1,\cdots,1)\\
0, & \mbox{ otherwise }.
\end{cases}
$$ 
\end{theorem}

 \section{Large Deviations for Small-time Sampling Distributions}
 
 Recall that 
 $$
  P_{n}^{\theta}(\eta)=\frac{n!}{\prod_{i=1}^{l(\eta)}\eta_{i}!\prod_{i=1}^{n}(i!)^{\alpha_{i}(\eta)}}\mathbb{E}p_{\eta}^{o}(D(Z_{t(\theta)})).
 $$
  By Theorem \ref{asymptotic_sampling}, we know the limiting sampling distribution is 
 \begin{equation}
 \frac{n!}{\prod_{i=1}^{l(\eta)}\eta_{i}!\prod_{i=1}^{n}(i!)^{\alpha_{i}}}p_{\eta}^{o}(D(m_{c})). \label{limit_sample}
 \end{equation}
  Let $D(m_{c})$ be $x=(x_{1},x_{2},\cdots)$. If $\sum_{i=1}^{\infty}x_{i}>0$, the randomness in the partition with distribution (\ref{limit_sample}) presents. Clearly, this randomness arises from sampling, not from allele frequency. When $\sum_{i=1}^{\infty}x_{i}=0$, the limiting sampling distribution is trivial. The randomness in the associated partition disappears. Therefore, when $\lim_{\theta\to\infty}t(\theta)=0$ and $\lim_{\theta}\theta t(\theta)=\infty$, one will end up with $m_{\infty}=\nu_{0}$, which is diffuse. Then $D(\nu_{0})$ is $x=(0,0,\cdots)$. The partition associated with random sample is always $(1,1,\cdots,1)$. In this period of evolution process, one should  only expect such trivial sample partition. Other partitions can be very rare. The degree of rareness can be estimated by large deviation principle(LDP in short). In this section, we are going to establish the large deviations for sampling distributions $\{P_{n}^{\theta}, \theta>0\}$ under the small-time scale $t(\theta)$ where $\lim_{\theta\to\infty}t(\theta)=0$ and $\lim_{\theta}\theta t(\theta)=\infty$.
   
 \begin{theorem}\label{partition_LDP}
When $\lim_{\theta\to\infty}t(\theta)=0$ and $\lim_{\theta}\theta t(\theta)=\infty$, as $\theta\rightarrow+\infty$, for any given integer $n\geq2$, $P^{\theta}_{n}$ satisfies an LDP with speed $\theta t(\theta)$ and rate function
$$
I(\eta)=\min\left\{\frac{n-l(\eta)}{k}, \frac{n-\alpha_{1}(\eta)}{2}\right\},
$$
where $k=\lim_{\theta\to\infty}\frac{\theta t(\theta)}{\log\theta}$.
\end{theorem}
\begin{rmk}
\item[1,] Since 
\begin{align*}
\frac{n-l(\eta)}{k}-\frac{n-\alpha_{1}(\eta)}{2}=&\frac{1}{k}\left[(n-l(\eta))-\frac{k}{2}(n-\alpha_{1}(\eta))\right]\\
=&\frac{1}{k}\left[(n-\alpha_{1}(\eta))-(l(\eta)-\alpha_{1}(\eta))-\frac{k}{2}(n-\alpha_{1}(\eta))\right]\\
=&\frac{1}{k}(l(\eta)-\alpha_{1}(\eta))\left[\frac{n-\alpha_{1}(\eta)}{l(\eta)-\alpha_{1}(\eta)}\frac{2-k}{2}-1\right],
\end{align*}
the value of $I(\eta)$ depends solely on $\frac{n-\alpha_{1}(\eta)}{l(\eta)-\alpha_{1}(\eta)}$, which is an important index of species diversity in population genetics $($cf. \cite{MR0368829}$)$.
\item[2,] The rate function indicates a competition of two degrees of rareness.  One is of order $\exp\{-\theta t(\theta)\frac{(n-l(\eta))}{k}\}$, which is also equivalent to $\exp\{-\log(\theta)(n-l(\eta))\}$; another is of order $\exp\{-\theta t(\theta)\frac{(n-\alpha_{1}(\eta))}{2}\}$. Thus, the competition of these two effects will result in phase transition. Note that the LDP for Ewens sampling distribution $(\ref{Ewens_Sampling_Formula})$ in \cite{MR2663265} has speed $\log\theta$ and rate function $n-l(\eta)$ as well. So when $k\geq2$, the random sample has the same degree of rareness as equilibrium case. When $1\leq k<2$, new features start to emerge, and the two degrees of rareness are entangled. When $0\leq k<1$, the new feature dominates. 
\item[3,] Making use of the transition density of the infinite dimensional diffusion associated with the two-parameter Poisson-Dirichlet distribution in \cite{MR2737405}, we can obtain the similar result.
\end{rmk}

The proof of Theorem \ref{partition_LDP} depends heavily on the following lemma.
 \begin{lemma}\label{Inner_estimation}
  As $\theta\rightarrow+\infty$, for any given partition $\eta\in\mathcal{J}^{o}$, we have the following asymptotic estimations:
  \begin{equation}
  \langle\varphi_{\eta},1\rangle_{\theta}\sim(\eta_{1}-1)!\cdots(\eta_{l}-1)!\frac{1}{\theta^{|\eta|-l(\eta)}}. \label{first_part}
 \end{equation}
 For a given partition $\xi\in\mathcal{J}^{o}$, we have $\forall \eta\geq\xi,\eta\in\mathcal{J}^{o}$
 \begin{align*}
 \langle\varphi_{\eta},\psi_{\xi}^{\theta}\rangle_{\theta}\sim&\left[\sum_{i=1}^{l(\eta)}\sum_{j=1}^{p(\xi)}\frac{(\eta_{i}+\xi_{j}-1)!}{(\eta_{i}-1)!(\xi_{j}-1)!}-|\eta||\xi|\right]
 (\eta_{1}-1)!\cdots(\eta_{l(\eta)}-1)!\\
 &(\xi_{1}-1)!\cdots(\xi_{l(\xi)}-1)!\frac{1}{\theta^{|\eta|-l(\eta)+|\xi|-l(\xi)+1}}.
 \end{align*}
  \end{lemma}
  
 [\textbf{Proof of Theorem \ref{partition_LDP}}]:  
 \begin{proof}
According to Theorem \ref{asymptotic_sampling}
$$
P_{n}^{\theta}(\eta)\rightarrow\delta_{(1,1,\cdots,1)}(\eta).
$$
Then $\log P_{n}^{\theta}(1,1,\cdots,1)\rightarrow0$, as $\theta\rightarrow+\infty$. So we only need to consider 
$\eta\neq(1,1,\cdots,1)$.
Suppose that $|\eta|=n$. By Theorem \ref{sampling_distribution_theorem}, we have
$$
P_{n}^{\theta}(\eta)=\frac{n!}{\eta_{1}!\cdots\eta_{l(\eta)}!\alpha_{1}(\eta)!\cdots\alpha_{n}(\eta)!}\mathbb{E}p^{o}_{\eta}(D(Z_{t(\theta)})),
$$
where 
$$
p_{\eta}^{o}(x)=\sum_{d=1}^{l(\eta)}(-1)^{l(\eta)-d}\sum_{\beta\in\pi(l(\eta),d)}(|\beta_{1}|-1)!\cdots (|\beta_{d}|-1)! \varphi_{\sum_{i\in\beta_{1}}\eta_{i}}\cdots\varphi_{\sum_{i\in\beta_{d}}\eta_{i}}.
$$
 Let $\gamma$ be the decreasing arrangement of $(\sum_{i\in\beta_{1}}\eta_{i},\cdots,\sum_{i\in\beta_{d}}\eta_{i})$. We then define $\xi^{\beta}$ to be the projection of $\gamma$ onto $\mathcal{J}^{o}$. Note that we always have $|\xi^{\beta}|\geq n-\alpha_{1}(\eta)$ and
 \begin{equation}
|\xi|>|\xi|-l(\xi)=n-d\geq1,\quad \forall \xi\in\mathcal{J}^{o}\label{index}
 \end{equation}
  By the transition density (\ref{density}) of $X_{t}=D(Z_{t})$ , we know
$$
\mathbb{E}p^{o}_{\eta}(D(Z_{t(\theta)}))=\int_{\bar{\triangledown}_{\infty}}p_{\eta}^{o}(y)\mbox{PD}(\theta)(dy)+\sum_{m=2}^{\infty}e^{-\lambda_{m}t(\theta)}\int_{\bar{\triangledown}_{\infty}}p_{\eta}^{o}(y)Q_{m}(x,y)\mbox{PD}(\theta)(dy).
$$
We define 
$$
A_{1}=\int_{\bar{\triangledown}_{\infty}}p_{\eta}^{o}(y)\mbox{PD}(\theta)(dy)
$$
and 
$$
A_{2}=\sum_{m=2}^{\infty}e^{-\lambda_{m}t(\theta)}\int_{\bar{\triangledown}_{\infty}}p_{\eta}^{o}(y)Q_{m}(x,y)\mbox{PD}(\theta).
$$
By Ewens sampling distribution (\ref{Ewens_Sampling_Formula}), one easily show that 
\begin{equation}
\begin{split}
A_{1}=&(\eta_{1}-1)!\cdots(\eta_{l(\eta)}-1)!\frac{\theta^{l(\eta)}}{\theta_{(n)}}\\
\sim&(\eta_{1}-1)!\cdots(\eta_{l(\eta)}-1)!\exp\{-(\log\theta) (n-l(\eta))\} \label{A1}
\end{split}
\end{equation}
Now it remains to calculate $A_{2}$. If we can pinpoint the leading term in $A_{1}+A_{2}$, then the LDP for $P_{n}^{\theta}$ can be followed readily. By the definition of $Q_{m}(x,y)$ and the orthogonalization process with $\varphi_{\eta}$, we can conclude that the infinite sum in $A_{2}$ is essentially a finite sum for there is a cutoff at the $n$th term. So it is possible to determine the leading term in $A_{1}+A_{2}$. 

Substituting the explicit expression of $p_{\eta}^{o}(y)$ into $A_{2}$,  
we have
\begin{align*}
A_{2}=\sum_{d=1}^{l(\eta)}(-1)^{l(\eta)-d}&\sum_{\beta\in\pi(l(\eta),d)}\prod_{i=1}^{d}(|\beta_{i}|-1)!\sum_{m=2}^{\infty}e^{-\lambda_{m}t(\theta)}\\
&\int_{\bar{\triangledown}_{\infty}}\varphi_{\xi^{\beta}}(y)Q_{m}(x,y)\mbox{PD}(\theta)(dy)
\end{align*} 
By the Gram-Schmidt orthogonalization process with $\varphi_{\xi^{\beta}}$, we have 
\begin{equation}
\varphi_{\xi^{\beta}}=\psi_{\xi^{\beta}}^{\theta}+\sum_{\delta<\xi^{\beta},\delta\in\mathcal{J}}\langle\varphi_{\xi^{\beta}},\psi_{\delta}^{\theta}\rangle_{\theta}\psi_{\delta}^{\theta}.\label{GS}
\end{equation}
Due to Lemma \ref{orthonormal}, replacing $\varphi_{\xi^{\beta}}$ in $A_{2}$ by the right hand side of (\ref{GS}) yields 
\begin{align*}
&\sum_{d=1}^{l(\eta)}(-1)^{l(\eta)-d}\sum_{\beta\in\pi(l(\eta),d)}\prod_{i=1}^{d}(|\beta_{i}|-1)!
\sum_{m=2}^{\infty}e^{-\lambda_{m}t(\theta)}\\
&\left(\psi_{\xi^{\beta}}^{\theta}(x)\delta_{m,|\xi^{\beta}|}+\sum_{\delta<\xi^{\beta},\delta\in\mathcal{J}^{o}}\langle\varphi_{\xi^{\beta}},\psi_{\delta}^{\theta}\rangle_{\theta}\psi_{\delta}^{\theta}(x)\delta_{m,|\delta|}\right)\\
=&\sum_{d=1}^{l(\eta)}(-1)^{l(\eta)-d}\sum_{\beta\in\pi(l(\eta),d)}\prod_{i=1}^{d}(|\beta_{i}|-1)!\\
&\left[\psi_{\xi^{\beta}}^{\theta}(x)e^{-\lambda_{|\xi^{\beta}|}t(\theta)}+\sum_{\delta<\xi^{\beta},\delta\in\mathcal{J}^{o}}\langle\varphi_{\xi^{\beta}},\psi_{\delta}^{\theta}\rangle_{\theta}\psi_{\delta}^{\theta}(x)e^{-\lambda_{|\delta|}t(\theta)}\right].
\end{align*}
By (\ref{GS}), we rewrite $\psi_{\xi^{\beta}}^{\theta}$ as $\varphi_{\xi^{\beta}}-\sum_{\delta<\xi^{\beta},\delta\in\mathcal{J}}\langle\varphi_{\xi^{\beta}},\psi_{\delta}^{\theta}\rangle_{\theta}\psi_{\delta}^{\theta}$ and substitute it into $A_{2}$. Then
$$
A_{2}=M_{1}+M_{2}+M_{3},
$$
where 
$$
M_{1}=\sum_{d=1}^{l(\eta)}(-1)^{l(\eta)-d}\sum_{\beta\in\pi(l(\eta),d)}\prod_{i=1}^{d}(|\beta_{i}|-1)!
\varphi_{\xi^{\beta}}(x)e^{-\lambda_{|\xi^{\beta}|}t(\theta)},
$$
\begin{align*}
M_{2}=&\sum_{d=1}^{l(\eta)}(-1)^{l(\eta)-d}\sum_{\beta\in\pi(l(\eta),d)}\prod_{i=1}^{d}(|\beta_{i}|-1)!\\
&\sum_{\delta<\xi^{\beta},\delta\in\mathcal{J}^{o}}\langle\varphi_{\xi^{\beta}},\psi_{\delta}^{\theta}\rangle_{\theta}\psi_{\delta}^{\theta}(x)\left(e^{-\lambda_{|\delta|}t(\theta)}-e^{-\lambda_{|\xi^{\beta}|}t(\theta)}\right)
\end{align*}
and
$$
M_{3}=-\sum_{d=1}^{l(\eta)}(-1)^{l(\eta)-d}\sum_{\beta\in\pi(l(\eta),d)}\prod_{i=1}^{d}(|\beta_{i}|-1)!\langle\varphi_{\xi^{\beta}},1\rangle_{\theta}e^{-\lambda_{|\xi^{\beta}|t(\theta)}}.
$$
 Recall that $|\xi^{\beta}|\geq n-\alpha_{1}(\eta)$,  the equality holds when 
 $$
 \beta=\beta^{'}\cup\{l(\eta)-\alpha_{1}(\eta)+1\}\cup\cdots\cup\{l(\eta)\},
 $$
 where $\beta^{'}\in\pi(l(\eta)-\alpha_{1}(\eta),d^{'}), 1\leq d^{'}\leq l(\eta)-\alpha_{1}(\eta)$.
Thus, we have
\begin{align*}
M_{1}\sim&\Big[\sum_{d=1}^{l(\eta)-\alpha_{1}(\eta)}(-1)^{l(\eta)-(d+\alpha_{1}(\eta))}\sum_{\beta\in\pi(l(\eta)-\alpha_{1}(\eta),d)}\prod_{i=1}^{d}(|\beta_{i}|-1)!\varphi_{\xi^{\beta}}(x)\Big]\\
&e^{-\lambda_{n-\alpha_{1}(\eta)}t(\theta)}\\
=&p_{\eta_{1},\cdots,\eta_{l-\alpha_{1}(\eta)}}^{o}(x)e^{-\lambda_{n-\alpha_{1}(\eta)}t(\theta)},\quad\mbox{ due to proposition \ref{partition_representation}.}\\
\sim &p_{\eta_{1},\cdots,\eta_{l-\alpha_{1}(\eta)}}^{o}(x)e^{-\frac{n-\alpha_{1}(\eta)}{2}\theta t(\theta)}
\end{align*}

By Lemma \ref{Inner_estimation}, we know $\forall\delta<\xi^{\beta}(\delta\in\mathcal{J}^{o})$,
\begin{align*}
\langle\varphi_{\xi^{\beta}},\psi_{\delta}^{\theta}\rangle_{\theta}\sim&\Bigg[\sum_{i=1}^{l(\xi^{\beta})}\sum_{j=1}^{l(\delta)}\frac{(\xi^{\beta}_{i}+\delta_{j}-1)!}{(\xi^{\beta}_{i}-1)!(\delta_{j}-1)!}-|\xi^{\beta}||\delta|\Bigg]
(\xi_{1}^{\beta}-1)!\cdots(\xi_{d}^{\beta}-1)!\\
&(\delta_{1}-1)!\cdots(\delta_{l(\delta)}-1)!\frac{1}{\theta^{|\xi^{\beta}|-l(\xi^{\beta})+|\delta|-l(\delta)+1}}
\end{align*}
Recall that for $\xi^{\beta},\delta\in\mathcal{J}^{o}$ we have (\ref{index}). Then 
$$
|\xi^{\beta}|-l(\xi^{\beta})+|\delta|-l(\delta)+1=n-d+|\delta|-l(\delta)+1\geq n-d+2.
$$
Moreover, for $1\leq d\leq l(\eta)$, $|\xi^{\beta}|-l(\xi^{\beta})+|\delta|-l(\delta)+1\geq n-l(\eta)+2$. The equality holds if and only if $d=l(\eta)$ and $\delta=(2)$. So if $\xi^{\beta}=(2)$, then $M_{2}=0$; if $\xi^{\beta}>(2),$ then 
$$
M_{2}\sim \Big[\sum_{d=1}^{l(\eta)-\alpha_{1}(\eta)}\frac{(\eta_{i}+1)!}{(\eta_{i}-1)!}-2(n-\alpha_{1}(\eta))\Big]\prod_{i=1}^{l(\eta)-\alpha_{1}(\eta)}(\eta_{i}-1)!\varphi_{2}\frac{1}{\theta^{n-l(\eta)+2}}e^{-(\theta+1)t(\theta)}
$$
and
 \begin{align*}
M_{3}\sim& -(\eta_{1}-1)!\cdots(\eta_{l(\eta)-\alpha_{1}(\eta)}-1)!\frac{1}{\theta^{n-l(\eta)}}e^{-\frac{n-\alpha_{1}(\eta)}{2}\theta t(\theta)}
\end{align*}
Obviously,  $\lim_{\theta\to\infty}\frac{M_{2}}{A_{1}}=\lim_{\theta\to\infty}\frac{M_{3}}{A_{1}}=0$. Then 
$$
A_{1}+A_{2}=A_{1}\left(1+\frac{M_{2}}{A_{1}}+\frac{M_{3}}{A_{1}}\right)+M_{1}.
$$
Note that
$$
\frac{1}{\theta t(\theta)}\log P_{n}^{\theta}(\eta)\sim \frac{1}{\theta t(\theta)}\log (A_{1}+A_{2}).
$$
 By Lemma 1.2.15 in \cite{MR2571413}, we have
\begin{align*}
&\limsup_{\theta\to\infty}\frac{1}{\theta t(\theta)}\log (A_{1}+A_{2})\\
\leq &\max\left\{\limsup_{\theta\to\infty}\frac{1}{\theta t(\theta)}\log \left[A_{1}\left(1+\frac{M_{2}}{A_{1}}+\frac{M_{3}}{A_{1}}\right)\right],\limsup_{\theta\to\infty}\frac{1}{\theta t(\theta)}\log M_{1}\right\}\\
\leq& \max\left\{\limsup_{\theta\to\infty}\frac{1}{\theta t(\theta)}\log A_{1},\limsup_{\theta\to\infty}\frac{1}{\theta t(\theta)}\log M_{1}\right\}
\end{align*}
Moreover, for large enough $\theta$, both $A_{1}\left(1+\frac{M_{2}}{A_{1}}+\frac{M_{3}}{A_{1}}\right)$ and $M_{1}$ are positive. So we also have
\begin{align*}
&\liminf_{\theta\to\infty}\frac{1}{\theta t(\theta)}\log (A_{1}+A_{2})\\
\geq &\max\left\{\liminf_{\theta\to\infty}\frac{1}{\theta t(\theta)}\log \left [A_{1}\left(1+\frac{M_{2}}{A_{1}}+\frac{M_{3}}{A_{1}}\right)\right],\liminf_{\theta\to\infty}\frac{1}{\theta t(\theta)}\log M_{1}\right\}\\
\geq & \max\left\{\liminf_{\theta\to\infty}\frac{1}{\theta t(\theta)}\log A_{1},\liminf_{\theta\to\infty}\frac{1}{\theta t(\theta)}\log M_{1}\right\}\end{align*}
If $\lim_{\theta\to\infty}\frac{\theta t(\theta)}{\log\theta}=k$, we can easily show that 
$$
\lim_{\theta\to\infty}\frac{1}{\theta t(\theta)}\log A_{1}=-\frac{n-l(\eta)}{k}
$$
and
$$
\lim_{\theta\to\infty}\frac{1}{\theta t(\theta)}\log M_{1}=-\frac{n-\alpha_{1}(\eta)}{2}.
$$
Therefore,
$$
\lim_{\theta\to\infty}\frac{1}{\theta t(\theta)}\log P_{n}^{\theta}(\eta)=\lim_{\theta\to\infty}\frac{1}{\theta t(\theta)}\log (A_{1}+A_{2})=-\min\left\{\frac{n-l(\eta)}{k},\frac{n-\alpha_{1}(\eta)}{2}\right\}.
$$

The LDP for $P_{n}^{\theta}$ is established. 
\end{proof}

 [\textbf{Proof of Lemma \ref{Inner_estimation}}] :

 \begin{proof}
 For partition $\eta\in\mathcal{J}^{o}$ and $\beta\in\pi(l(\eta),d)$, we define $\eta^{\beta}$ to be the decreasing arrangement of $(\sum_{i\in\beta_{1}}\eta_{i},\cdots,\sum_{i\in\beta_{d}}\eta_{i}).$ Since $\varphi_{\eta}=\sum_{d=1}^{l(\eta)}\sum_{\beta\in\pi(l(\eta),d)}p^{o}_{\eta^{\beta}}$,we have 
 \begin{align*}
 \langle\varphi_{\eta},1\rangle_{\theta}=\sum_{d=1}^{l(\eta)}\sum_{\beta\in\pi(l(\eta),d)}\int_{\bar{\triangledown}_{\infty}}p^{o}_{\eta^{\beta}}d\mbox{PD}(\theta).
\end{align*}
By Ewens sampling distribution (\ref{Ewens_Sampling_Formula}), we have
\begin{align}
\langle\varphi_{\eta},1\rangle_{\theta}=\sum_{d=1}^{l(\eta)}\sum_{\beta\in\pi(l(\eta),d)}
(\sum_{i\in\beta_{1}}\eta_{i}-1)!\cdots(\sum_{i\in\beta_{d}}\eta_{i}-1)!\frac{\theta^{d}}{\theta_{(|\eta|)}},\label{Leading_term}.
\end{align} 
 Therefore, the leading term in  (\ref{Leading_term}) is the term $\int_{\bar{\triangledown}_{\infty}}p^{o}_{\eta^{\beta}}d\mbox{PD}(\theta)$ where  $\beta=\{1\}\cup\{2\}\cup\cdots\{l(\eta)\}$. Then 
 $$
 \langle\varphi_{\eta},1\rangle_{\theta}\sim(\eta_{1}-1)!\cdots(\eta_{l(\eta)}-1)!\frac{1}{\theta^{|\eta|-l(\eta)}}. 
 $$

 The first statement is thus proved.   Now we can use mathematical induction with respect to partition $\xi (\xi\in\mathcal{J}^{o})$ to show the second statement. We will complete this in three steps.\\
 \begin{itemize}
\item[1] \textbf{The initial case:} $\langle\varphi_{\eta},\psi_{2}^{\theta}\rangle_{\theta},\forall\eta\in\mathcal{J}^{o}$ and $\eta\geq(2)$\\
Since $\psi_{2}^{\theta}=\varphi_{2}-\frac{1}{1+\theta}$, we know
  \begin{equation}
  \langle\varphi_{\eta},\psi_{2}^{\theta}\rangle_{\theta}=\langle\varphi_{\eta},\varphi_{2}\rangle_{\theta}-\frac{1}{1+\theta}\langle\varphi_{\eta},1\rangle_{\theta}\label{Check_2}
 \end{equation} 
  Define $\tau$ to be the decreasing arrangement of $(\eta_{1},\cdots,\eta_{l(\eta)},2)$. Then
 $$\langle\varphi_{\eta},\varphi_{2}\rangle_{\theta}=\langle\varphi_{\tau},1\rangle_{\theta}.$$ 
 Using Ewens sampling formula (\ref{Ewens_Sampling_Formula}), we have
 \begin{equation}
 \langle\varphi_{\eta},\varphi_{2}\rangle_{\theta}=\sum_{d=1}^{l(\eta)+1}\sum_{\beta\in\pi(l(\eta)+1,d)}(\sum_{i\in\beta_{1}}\tau_{i}-1)!\cdots(\sum_{i\in\beta_{d}}\tau_{i}-1)!\frac{\theta^{d}}{\theta_{(|\eta|+2)}}. \label{first_check_2}
  \end{equation}
Substituting (\ref{first_check_2}) and (\ref{Leading_term}) into (\ref{Check_2}) yields $\langle\varphi_{\eta},\psi_{2}^{\theta}\rangle_{\theta}=J_{1}-J_{2}$, where
$$
J_{1}=\sum_{d=1}^{l(\eta)+1}\sum_{\beta\in\pi(l(\eta)+1,d)}(\sum_{i\in\beta_{1}}\tau_{i}-1)!\cdots(\sum_{i\in\beta_{d}}\tau_{i}-1)!\frac{\theta^{d}}{\theta_{(|\eta|+2)}}
$$
and
$$
J_{2}=\sum_{d=1}^{l(\eta)}\sum_{\beta\in\pi(l(\eta),d)}
(\sum_{i\in\beta_{1}}\eta_{i}-1)!\cdots(\sum_{i\in\beta_{d}}\eta_{i}-1)!\frac{\theta^{d}}{\theta_{(|\eta|)}}\frac{1}{1+\theta}.
$$
Then we split $J_{1}$ and $J_{2}$ into two parts as follows:
$$
J_{1}=R_{1}^{1}+R_{2}^{1},\quad J_{2}=R_{1}^{2}+R_{2}^{2},
$$
where 
\begin{align*}
R_{1}^{1}=&\sum_{d=l(\eta)}^{l(\eta)+1}\sum_{\beta\in\pi(l(\eta)+1,d)}(\sum_{i\in\beta_{1}}\tau_{i}-1)!\cdots(\sum_{i\in\beta_{d}}\tau_{i}-1)!\frac{\theta^{d}}{\theta_{(|\eta|+2)}}\\
R_{2}^{1}=&\sum_{d=1}^{l(
\eta)-1}\sum_{\beta\in\pi(l(\eta)+1,d)}(\sum_{i\in\beta_{1}}\tau_{i}-1)!\cdots(\sum_{i\in\beta_{d}}\tau_{i}-1)!\frac{\theta^{d}}{\theta_{(|\eta|+2)}}\\
R_{1}^{2}=&\sum_{d=l(\eta)-1}^{l(\eta)}\sum_{\beta\in\pi(l(\eta),d)}
(\sum_{i\in\beta_{1}}\eta_{i}-1)!\cdots(\sum_{i\in\beta_{d}}\eta_{i}-1)!\frac{\theta^{d}}{\theta_{(|\eta|)}}\frac{1}{1+\theta}\\
R_{2}^{2}=&\sum_{d=1}^{l(\eta)-2}\sum_{\beta\in\pi(l(\eta),d)}
(\sum_{i\in\beta_{1}}\eta_{i}-1)!\cdots(\sum_{i\in\beta_{d}}\eta_{i}-1)!\frac{\theta^{d}}{\theta_{(|\eta|)}}\frac{1}{1+\theta}.\end{align*}
Thus $R_{2}^{1}$ and $R_{2}^{2}$ are at least of order $\frac{1}{\theta^{|\eta|+3-l(\eta)}}$. We claim that 
\begin{equation}
R_{1}^{1}-R_{1}^{2}\sim (\eta_{1}-1)!\cdots(\eta_{l(\eta)}-1)!\left(\sum_{u=1}^{l(\eta)}\eta_{u}(\eta_{u}+1)-2|\eta|\right)\frac{1}{\theta^{|\eta|+2-l(\eta)}}.\label{R1R2}
\end{equation} 
Therefore, for the initial case $\langle\varphi_{\eta},\psi_{2}^{\theta}\rangle_{\theta}$, we have
$$
 \langle\varphi_{\eta},\psi_{2}^{\theta}\rangle_{\theta}\sim (\eta_{1}-1)!\cdots(\eta_{l(\eta)}-1)!\left(\sum_{u=1}^{l(\eta)}\eta_{u}(\eta_{u}+1)-2|\eta|\right)\frac{1}{\theta^{|\eta|+2-l(\eta)}}. 
 $$
 
 Now we are going to show the claim (\ref{R1R2}).  To this end, we group $R_{1}^1-R_{2}^{1}$ as $K_{1}+K_{2}+K_{3}+K_{4}+K_{5}$, where
\begin{align*}
K_{1}=&(\eta_{1}-1)!\cdots(\eta_{l(\eta)}-1)!\frac{\theta^{l(\eta)+1}}{\theta_{(|\eta|+2)}}\\
K_{2}=&
 (\eta_{1}-1)!\cdots(\eta_{l(\eta)}-1)!\sum_{u=1}^{l(\eta)}(\eta_{u}+1)\eta_{u}\frac{\theta^{l(\eta)}}{\theta_{(|\eta|+2)}}\label{useful1}\\
K_{3}=& (\eta_{1}-1)!\cdots(\eta_{l(\eta)}-1)!\sum_{1\leq u<v\leq l(\eta)}\frac{(\eta_{u}+\eta_{v}-1)!}{(\eta_{u}-1)!(\eta_{v}-1)!}\frac{\theta^{l(\eta)}}{\theta_{(|\eta|+2)}}\\ 
K_{4}=&- (\eta_{1}-1)!\cdots(\eta_{l(\eta)}-1)!\sum_{1\leq u<v\leq l(\eta)}\frac{(\eta_{u}+\eta_{v}-1)!}{(\eta_{u}-1)!(\eta_{v}-1)!}\frac{\theta^{l(\eta)-1}}{\theta_{(|\eta|)}(\theta+1)}\\ 
K_{5}=&-(\eta_{1}-1)!\cdots(\eta_{l(\eta)}-1)!\frac{\theta^{l(\eta)}}{\theta_{(|\eta|)}}\frac{1}{\theta+1}.
\end{align*}
By arithmetic calculation, we can show that
\begin{align*}
K_{3}+K_{4}=-(\eta_{1}-1)!\cdots(\eta_{l(\eta)}-1)!&\sum_{1\leq u<v\leq l(\eta)}\frac{(\eta_{u}+\eta_{v}-1)!}{(\eta_{u}-1)!(\eta_{v}-1)!}\frac{\theta^{l(\eta)}}{\theta_{(|\eta|+2)}}\\
&\frac{2\theta|\eta|+|\eta|(|\eta|+1)}{\theta(\theta+1)},
\end{align*}
and 
\begin{align*}
&K_{1}+K_{2}+K_{5}\\
=&(\eta_{1}-1)!\cdots(\eta_{l(\eta)}-1)!\frac{\theta^{l(\eta)}}{\theta_{(|\eta|+2)}}\left[
\sum_{u=1}^{l(\eta)}\eta_{u}(\eta_{u}+1)+\theta-\frac{(|\eta|+\theta)(|\eta|+\theta+1)}{\theta+1}\right]\\
=&(\eta_{1}-1)!\cdots(\eta_{l(\eta)}-1)!\frac{\theta^{l(\eta)}}{\theta_{(|\eta|+2)}}\left[\sum_{u=1}^{l(\eta)}\eta_{u}(\eta_{u}+1)-\frac{2\theta+1}{\theta+1}|\eta|-\frac{|\eta|^{2}}{1+\theta}\right]\\
\sim&  (\eta_{1}-1)!\cdots(\eta_{l(\eta)}-1)!\left(\sum_{u=1}^{l(\eta)}\eta_{u}(\eta_{u}+1)-2|\eta|\right)\frac{1}{\theta^{|\eta|+2-l(\eta)}}.
\end{align*}
 Apparently, $K_{3}+K_{4}$ is of order $\frac{1}{\theta^{|\eta|-l(\eta)+3}}$, so $\lim_{\theta\to\infty}\frac{K_{3}+K_{4}}{K_{1}+K_{2}+K_{5}}=0$. Then 
 \begin{align*}
 R_{1}^{1}-R_{1}^{2}=&K_{1}+K_{2}+K_{3}+K_{4}+K_{5}\\
 =& (K_{1}+K_{2}+K_{5})\left(1+\frac{K_{3}+K_{4}}{K_{1}+K_{2}+K_{5}}\right)\\
 \sim& K_{1}+K_{2}+K_{5}\\
 \sim &(\eta_{1}-1)!\cdots(\eta_{l(\eta)}-1)!\left(\sum_{u=1}^{l(\eta)}\eta_{u}(\eta_{u}+1)-2|\eta|\right)\frac{1}{\theta^{|\eta|+2-l(\eta)}}.
 \end{align*}
The claim (\ref{R1R2}) is thus shown.\\
\item[2] \textbf{Induction Conjecture}: \\
For a given partition $\xi\in\mathcal{J}^{o},\xi>(2)$, we assume that $\forall\delta<\xi (\delta\in\mathcal{J}^{o})$ the following is true for all $\eta\geq\delta,\eta\in\mathcal{J}^{o}$ 
\begin{equation}
 \begin{split}
 \langle\varphi_{\eta},\psi_{\delta}^{\theta}\rangle_{\theta}\sim&\left[\sum_{i=1}^{l(\eta)}\sum_{j=1}^{l(\delta)}\frac{(\eta_{i}+\delta_{j}-1)!}{(\eta_{i}-1)!(\delta_{j}-1)!}-|\eta||\delta|\right]\\
 &(\eta_{1}-1)!\cdots(\eta_{l(\eta)}-1)!(\delta_{1}-1)!\cdots(\delta_{l(\delta)}-1)!
 \frac{1}{\theta^{|\eta|-l(\eta)+|\delta|-l(\delta)+1}}.\label{assumptionss}
 \end{split}
 \end{equation}
 \item[3] \textbf{Confirm the case  $\langle\varphi_{\eta},\psi_{\xi}^{\theta}\rangle_{\theta},\eta,\xi\in\mathcal{J}^{o},\eta\geq\xi$}.\\
 By (\ref{gram_schmidt}), we know
 $$
 \psi_{\xi}^{\theta}=\varphi_{\xi}-\sum_{\delta<\xi,\delta\in\mathcal{J}}\langle\varphi_{\xi},\psi_{\delta}^{\theta}\rangle_{\theta}\psi_{\delta}^{\theta}.
 $$
 Therefore,
$$
 \langle\varphi_{\eta},\psi_{\xi}^{\theta}\rangle_{\theta}
 =B_{1}+B_{2},
 $$
 where
 $$
B_{1}=\langle\varphi_{\eta},\varphi_{\xi}\rangle_{\theta}-\langle\varphi_{\xi},1\rangle_{\theta}\langle\varphi_{\eta},1\rangle_{\theta}
$$
and
$$
B_{2}=-\sum_{\delta<\xi,\delta\in\mathcal{J}^{o}}\langle\varphi_{\xi},\psi_{\delta}^{\theta}\rangle_{\theta}\langle\varphi_{\eta},\psi_{\delta}^{\theta}\rangle_{\theta}.
$$
By the above assumption (\ref{assumptionss}), for each partition $\delta<\xi,\delta\in\mathcal{J}^{o}$, the summand $\langle\varphi_{\xi},\psi_{\delta}^{\theta}\rangle_{\theta}\langle\varphi_{\eta},\psi_{\delta}^{\theta}\rangle_{\theta}$ in $B_{2}$ is of order 
$\frac{1}{\theta^{|\eta|-l(\eta)+|\xi|-l(\xi)+1}}\frac{1}{\theta^{2(|\delta|-l(\delta))+1}}$. Notice that $2(|\delta|-l(\delta))+1\geq 3$ for $\delta\in\mathcal{J}^{o}$. So the leading term in
$$
B_{2}=-\sum_{\delta<\xi,\delta\in\mathcal{J}^{o}}\langle\varphi_{\xi},\psi_{\delta}^{\theta}\rangle_{\theta}\langle\varphi_{\eta},\psi_{\delta}^{\theta}\rangle_{\theta}
$$
is of order $\frac{1}{\theta^{|\eta|+|\xi|-l(\eta)-l(\xi)+4}}$. 

Define $\omega$ to be the decreasing arrangement of $(\eta_{1},\cdots,\eta_{l(\eta)},\xi_{1},\cdots,\xi_{l(\xi)})$. Then $\langle\varphi_{\eta},\varphi_{\xi}\rangle_{\theta}=\langle\varphi_{\omega},1\rangle_{\theta}$. By Ewens sampling distribution (\ref{Ewens_Sampling_Formula}), we have
\begin{align*}
\langle\varphi_{\eta},\varphi_{\xi}\rangle_{\theta}=&\sum_{d=1}^{l(\eta)+l(\xi)}\sum_{\beta\in\pi(l(\eta)+l(\xi),d)}\int_{\bar{\triangledown}_{\infty}}p^{o}_{\omega}d\mbox{PD}(\theta)\\
=&\sum_{d=1}^{l(\eta)+l(\xi)}\sum_{\beta\in\pi(l(\eta)+l(\xi),d)}(\sum_{i\in\beta_{1}}\omega_{i}-1)!\cdots(\sum_{i\in\beta_{d}}\omega_{i}-1)!\frac{\theta^{d}}{\theta_{(|\eta|+|\xi|)}}.\\
\langle\varphi_{\eta},1\rangle_{\theta}=&\sum_{d=1}^{l(\eta)}\sum_{\beta\in\pi(l(\eta),d)}\int_{\bar{\triangledown}_{\infty}}p^{o}_{\eta}d\mbox{PD}(\theta)\\
=&\sum_{d=1}^{l(\eta)}\sum_{\beta\in\pi(l(\eta),d)}(\sum_{i\in\beta_{1}}\eta_{i}-1)!\cdots(\sum_{i\in\beta_{d}}\eta_{i}-1)!\frac{\theta^{d}}{\theta_{(|\eta|)}}\\
\langle\varphi_{\xi},1\rangle_{\theta}=&\sum_{d=1}^{l(\xi)}\sum_{\beta\in\pi(l(\xi),d)}\int_{\bar{\triangledown}_{\infty}}p^{o}_{\xi}d\mbox{PD}(\theta)\\
=&\sum_{d=1}^{l(\xi)}\sum_{\beta\in\pi(l(\xi),d)}(\sum_{i\in\beta_{1}}\xi_{i}-1)!\cdots(\sum_{i\in\beta_{d}}\xi_{i}-1)!\frac{\theta^{d}}{\theta_{(|\xi|)}}
\end{align*}
 
 We then split $B_{1}$ as $S_{1}+S_{2}+S_{3}+S_{4}+S_{5}+S_{6}$, where
 \begin{align*}
 S_{1}=&(\omega_{1}-1)!\cdots(\omega_{l(\eta)+l(\xi)}-1)!\frac{\theta^{l(\eta)+l(\xi)}}{\theta_{(|\eta|+|\xi|)}}\\
 S_{2}=&\left(\sum_{1\leq i<j\leq l(\eta)+l(\xi)}\frac{(\omega_{i}+\omega_{j}-1)!}{(\omega_{i}-1)!(\omega_{j}-1)!}\right)\prod_{i=1}^{l(\eta)+l(\xi)}(\omega_{i}-1)!\frac{\theta^{l(\eta)+l(\xi)-1}}{\theta_{(|\eta|+|\xi|)}}\\
 S_{3}=&-\prod_{i=1}^{l(\eta)}(\eta_{i}-1)!\prod_{j=1}^{l(\xi)}(\xi_{j}-1)!\frac{\theta^{l(\eta)+l(\xi)}}{\theta_{(|\eta|)}\theta_{(|\xi|)}}\\
 S_{4}=&-\left(\sum_{1\leq i< j\leq l(\eta)}\frac{(\eta_{i}+\eta_{j}-1)!}{(\eta_{i}-1)!(\eta_{j}-1)!} \right)\prod_{i=1}^{l(\eta)}(\eta_{i}-1)!\prod_{j=1}^{l(\xi)}(\xi_{j}-1)!\frac{\theta^{l(\eta)+l(\xi)-1}}{\theta_{(|\eta|)}\theta_{(|\xi|)}}\\
 S_{5}=&-\left(\sum_{1\leq i<j\leq l(\xi)}\frac{(\xi_{i}+\xi_{j}-1)!}{(\xi_{i}-1)!(\xi_{j}-1)!}\right)\prod_{i=1}^{l(\eta)}(\eta_{i}-1)!\prod_{j=1}^{l(\xi)}(\xi_{j}-1)!\frac{\theta^{l(\eta)+l(\xi)-1}}{\theta_{(|\eta|)}\theta_{(|\xi|)}}\\
 S_{6}=&\sum_{d=1}^{l(\eta)+l(\xi)-2}\sum_{\beta\in\pi(l(\eta)+l(\xi),d)}(\sum_{i\in\beta_{1}}\omega_{i}-1)!\cdots(\sum_{i\in\beta_{d}}\omega_{i}-1)!\frac{\theta^{d}}{\theta_{(|\eta|+|\xi|)}}\\
 &-\left[\sum_{d=1}^{l(\eta)-1}\sum_{\beta\in\pi(l(\eta),d)}(\sum_{i\in\beta_{1}}\eta_{i}-1)!\cdots(\sum_{i\in\beta_{d}}\eta_{i}-1)!\frac{\theta^{d}}{\theta_{(|\eta|)}}\right]\cdot \\
 &\left[\sum_{d=1}^{l(\xi)-1}\sum_{\beta\in\pi(l(\xi),d)}(\sum_{i\in\beta_{1}}\xi_{i}-1)!\cdots(\sum_{i\in\beta_{d}}\xi_{i}-1)!\frac{\theta^{d}}{\theta_{(|\xi|)}}\right]  
 \end{align*}
 
 It is easy to see that $S_{6}$ is of order $\frac{1}{\theta^{|\eta|+|\xi|-l(\eta)-l(\xi)+2}}$. We claim that 
 \begin{equation}
 \begin{split}
 &S_{1}+\cdots+S_{5}\sim \\
&\prod_{i=1}^{l(\eta)}(\eta_{i}-1)!\prod_{j=1}^{l(\xi)}(\xi_{j}-1)!\left[\sum_{i=1}^{l(\eta)}\sum_{j=1}^{l(\xi)}\frac{(\eta_{i}+\xi_{j}-1)!}{(\eta_{i}-1)!(\xi_{j}-1)!}
  -|\eta||\xi|\right]\frac{1}{\theta^{|\eta|-l(\eta)+|\xi|-l(\xi)+1}}. \label{S1S5}   \end{split}
  \end{equation}
   Then  $\lim_{\theta\to\infty}\frac{S_{6}}{S_{1}+\cdots+S_{5}}=0$. So by the claim (\ref{S1S5})
   \begin{align*}
   B_{1}=&S_{1}+\cdots+S_{6}=(S_{1}+\cdots+S_{5})\left(1+\frac{S_{6}}{S_{1}+\cdots+S_{5}}\right)\\
   \sim&S_{1}+\cdots+S_{5}\\
   \sim &\prod_{i=1}^{l(\eta)}(\eta_{i}-1)!\prod_{j=1}^{l(\xi)}(\xi_{j}-1)!\left[\sum_{i=1}^{l(\eta)}\sum_{j=1}^{l(\xi)}\frac{(\eta_{i}+\xi_{j}-1)!}{(\eta_{i}-1)!(\xi_{j}-1)!}
  -|\eta||\xi|\right]\\
  &\frac{1}{\theta^{|\eta|-l(\eta)+|\xi|-l(\xi)+1}}.   
   \end{align*}
   Then $\lim_{\theta\to\infty}\frac{B_{2}}{B_{1}}=0$. Therefore, 
  \begin{align*}
   &\langle\varphi_{\eta},\psi_{\xi}^{\theta}\rangle_{\theta}=B_{1}\left(1+\frac{B_{2}}{B_{1}}\right)\sim B_{1}\\
   \sim&\prod_{i=1}^{l(\eta)}(\eta_{i}-1)!\prod_{j=1}^{l(\xi)}(\xi_{j}-1)!\left[\sum_{i=1}^{l(\eta)}\sum_{j=1}^{l(\xi)}\frac{(\eta_{i}+\xi_{j}-1)!}{(\eta_{i}-1)!(\xi_{j}-1)!}
  -|\eta||\xi|\right]\frac{1}{\theta^{|\eta|-l(\eta)+|\xi|-l(\xi)+1}}.
  \end{align*}
  The the case $\langle\varphi_{\eta},\psi_{\xi}^{\theta}\rangle_{\theta}$ is thus confirmed. To complete the proof, we only need to show the claim (\ref{S1S5}). 
 
Notice that 
 $$
 (\omega_{1}-1)!\cdots(\omega_{l(\eta)+l(\xi)}-1)!=(\eta_{1}-1)!\cdots(\eta_{l(\eta)}-1)!(\xi_{1}-1)!\cdots(\xi_{l(\xi)}-1)!.
 $$
 Factoring out the term $\prod_{i=1}^{l(\eta)}(\eta_{i}-1)!\prod_{j=1}^{l(\xi)}(\xi_{j}-1)! \frac{\theta^{l(\eta)+l(\xi)-1}}{\theta_{(|\eta|+|\xi|)}}$ in $S_{1}+\cdots+S_{5}$ gives us  
 \begin{align*}
 &\prod_{i=1}^{l(\eta)}(\eta_{i}-1)!\prod_{j=1}^{l(\xi)}(\xi_{j}-1)! \frac{\theta^{l(\eta)+l(\xi)-1}}{\theta_{(|\eta|+|\xi|)}}\Bigg[\sum_{1\leq i<j\leq l(\eta)+l(\xi)}\frac{(\omega_{i}+\omega_{j}-1)!}{(\omega_{i}-1)!(\omega_{j}-1)!}\\
 &+\theta\left(1-\frac{\theta_{(|\eta|+|\xi|)}}{\theta_{(|\eta|)}\theta_{(|\xi|)}}\right)-\sum_{1\leq i<j\leq l(\eta)}\frac{(\eta_{i}+\eta_{j}-1)!}{(\eta_{i}-1)!(\eta_{j}-1)!}
 \frac{\theta_{(|\eta|+|\xi|)}}{\theta_{(|\eta|)}\theta_{(|\xi|)}}\\
 &-\sum_{1\leq i<j\leq l(\xi)}\frac{(\xi_{i}+\xi_{j}-1)!}{(\xi_{i}-1)!(\xi_{j}-1)!}
 \frac{\theta_{(|\eta|+|\xi|)}}{\theta_{(|\eta|)}\theta_{(|\xi|)}} \Bigg]\\
 \end{align*}
 By simplifying $1-\frac{\theta_{(|\eta|+|\xi|)}}{\theta_{(|\eta|)}\theta_{(|\xi|)}}$, we have
 \begin{align*}
1-\frac{\theta_{(|\eta|+|\xi|)}}{\theta_{(|\eta|)}\theta_{(|\xi|)}}=& \frac{\sum^{|\eta|}_{u=1}|\xi|^{u}\sum_{0\leq i_{1}<\cdots<i_{|\eta|-u}\leq |\eta|-1}(\theta+i_{1})\cdots(\theta+i_{|\eta|-u})}{\theta_{(|\eta|)}} \\
\sim & |\xi|\sum_{i=0}^{|\eta|-1}\frac{1}{\theta+i}.
 \end{align*}
 Then 
 \begin{align*}
 &S_{1}+\cdots+S_{5}\\
 \sim&\prod_{i=1}^{l(\eta)}(\eta_{i}-1)!\prod_{j=1}^{l(\xi)}(\xi_{j}-1)!\frac{1}{\theta^{|\eta|-l(\eta)+|\xi|-l(\xi)+1}}\Bigg[\sum_{1\leq i<j\leq l(\eta)+l(\xi)}\frac{(\omega_{i}+\omega_{j}-1)!}{(\omega_{i}-1)!(\omega_{j}-1)!}\\
 &-|\xi||\eta|-\sum_{1\leq i<j\leq l(\eta)}\frac{(\eta_{i}+\eta_{j}-1)!}{(\eta_{i}-1)!(\eta_{j}-1)!}-\sum_{1\leq i<j\leq l(\xi)}\frac{(\xi_{i}+\xi_{j}-1)!}{(\xi_{i}-1)!(\xi_{j}-1)!}\Bigg]  \\
 =&\prod_{i=1}^{l(\eta)}(\eta_{i}-1)!\prod_{j=1}^{l(\xi)}(\xi_{j}-1)!\Bigg[\sum_{i=1}^{l(\eta)}\sum_{j=1}^{l(\xi)}\frac{(\eta_{i}+\xi_{j}-1)!}{(\eta_{i}-1)!(\xi_{j}-1)!}
  -|\eta||\xi|\Bigg]\frac{1}{\theta^{|\eta|-l(\eta)+|\xi|-l(\xi)+1}}.  
  \end{align*}
 This claim (\ref{S1S5}) is thus proved.
  \end{itemize}
\end{proof}

\section{Appendix}

[\textbf{Proof of Proposition \ref{partition_representation}}]

 \begin{proof}
 We can prove this by mathematical induction on $l(\eta)$. When $l(\eta)=1$, it is trivial; when $l(\eta)=2$, 
 \begin{equation}
 p_{(\eta_{1},\eta_{2})}^{o}(x)=p^{o}_{(\eta_{1})}(x)p^{o}_{(\eta_{2})}(x)-p^{o}_{(\eta_{1}+\eta_{2})}(x)=\varphi_{\eta_{1}}\varphi_{\eta_{2}}-\varphi_{\eta_{1}+\eta_{2}}.\label{iteration_of_polynomial}
 \end{equation}
 Assume that, for $l(\eta)=l$, we have
  \begin{equation*}
 p_{\eta}^{o}=\sum_{d=1}^{l}(-1)^{l-d}\sum_{\beta\in\pi(l,d)}(|\beta_{1}|-1)!\cdots (|\beta_{d}|-1)! \varphi_{\sum_{i\in\beta_{1}}\eta_{i}}\cdots\varphi_{\sum_{i\in\beta_{d}}\eta_{i}}.
 \end{equation*} 
Then, for $l(\eta)=l+1$, 
 \begin{eqnarray}
 p_{(\eta_{1},\cdots,\eta_{l},\eta_{l+1})}^{o}&=&p_{(\eta_{1},\cdots,\eta_{l})}^{o}\varphi_{\eta_{l+1}}-\sum_{u=1}^{l}p_{(\eta_{1},\cdots,\eta_{u}+\eta_{l+1},\cdots,\eta_{l})}^{o}\nonumber\\
 &=&\sum_{d=1}^{l}(-1)^{l-d}\sum_{\beta\in\pi(l,d)}(|\beta_{1}|-1)!\cdots (|\beta_{d}|-1)! \varphi_{\sum_{i\in\beta_{1}}\eta_{i}}\cdots\varphi_{\sum_{i\in\beta_{d}}\eta_{i}}\varphi_{\eta_{l+1}}\nonumber\\
 &&-\sum_{u=1}^{l}\sum_{d=1}^{l}(-1)^{l-d}\sum_{\beta\in\pi(l,d)}(|\beta_{1}|-1)!\cdots (|\beta_{d}|-1)! \varphi_{\sum_{i\in\beta_{1}}\eta_{i}^{u}}\cdots\varphi_{\sum_{i\in\beta_{d}}\eta_{i}^{u}}, \label{second_term}
 \end{eqnarray}
 where 
 $$
 \eta_{i}^{u}=\begin{cases}
 \eta_{i}  &\mbox{ if } i\neq u\\
 \eta_{i}+\eta_{l+1}  &\mbox{ if } i=u.
 \end{cases}
 $$
 By switching the order of summation in (\ref{second_term}), we have
 \begin{align*}
 &\sum_{u=1}^{l}\sum_{d=1}^{l}(-1)^{l-d}\sum_{\beta\in\pi(l,d)}(|\beta_{1}|-1)!\cdots (|\beta_{d}|-1)! \varphi_{\sum_{i\in\beta_{1}}\eta_{i}^{u}}\cdots\varphi_{\sum_{i\in\beta_{d}}\eta_{i}^{u}}\\
  =&\sum_{d=1}^{l}(-1)^{l-d}\sum_{\beta\in\pi(l,d)}(|\beta_{1}|-1)!\cdots (|\beta_{d}|-1)! \sum_{u=1}^{l}\varphi_{\sum_{i\in\beta_{1}}\eta_{i}^{u}}\cdots\varphi_{\sum_{i\in\beta_{d}}\eta_{i}^{u}},
 \end{align*}
 where, as a matter of fact, 
 \begin{align*}
 \sum_{u=1}^{l}\varphi_{\sum_{i\in\beta_{1}}\eta_{i}^{u}}\cdots\varphi_{\sum_{i\in\beta_{d}}\eta_{i}^{u}}=&\sum_{v=1}^{d}\sum_{u\in\beta_{v}}\varphi_{\sum_{i\in\beta_{1}}\eta_{i}^{u}}\cdots\varphi_{\sum_{i\in\beta_{d}}\eta_{i}^{u}}\\
 =&\sum_{v=1}^{d}|\beta_{v}| \varphi_{\sum_{i\in\beta_{1}}\eta_{i}}\cdots\varphi_{\sum_{i\in\beta_{v}}\eta_{i}+\eta_{l+1}}\cdots\varphi_{\sum_{i\in\beta_{d}}\eta_{i}}.
 \end{align*}
 Therefore, (\ref{second_term}) becomes
 \begin{align*}
 &\sum_{d=1}^{l}(-1)^{l-d}\sum_{v=1}^{d}\sum_{\beta\in\pi(l,d)}(|\beta_{1}|-1)!\cdots(|\beta_{v-1}|-1)! |\beta_{v}|!(|\beta_{v+1}|-1)!\cdots(|\beta_{d}|-1)! \\
 &\varphi_{\sum_{i\in\beta_{1}}\eta_{i}}\cdots\varphi_{\sum_{i\in\beta_{v}}\eta_{i}+\eta_{l+1}}\cdots\varphi_{\sum_{i\in\beta_{d}}\eta_{i}} \\
 =&\sum_{d=1}^{l}(-1)^{l-d}\sum_{v=1}^{d} \sum_{\beta^{v}\in\pi(l+1,d)}(|\beta_{1}^{v}|-1)!\cdots(|\beta_{d}^{v}|-1)!\varphi_{\sum_{i\in\beta_{1}^{v}}\eta_{i}}\cdots\varphi_{\sum_{i\in\beta_{d}^{v}}\eta_{i}}. 
 \end{align*}
 Here $\beta^{v}=\{\beta_{1}^{v},\cdots,\beta_{d}^{v}\}$,
 $$
 \beta_{i}^{v}=\begin{cases}
 \beta_{i} & \mbox{ if } i\neq v\\
 \beta_{i}\cup\{l+1\} &\mbox{ if }i=v,
 \end{cases}
 $$
 and $\beta=\{\beta_{1},\cdots,\beta_{d}\}\in\pi(l,d)$. Thus, $p_{(\eta_{1},\cdots,\eta_{l},\eta_{l+1})}^{o}$ is
 \begin{align}
&p_{(\eta_{1},\cdots,\eta_{l},\eta_{l+1})}^{o}\nonumber\\
 =&\sum_{d=2}^{l+1}(-1)^{l+1-d}\sum_{\beta\in\pi(l+1,d),\beta_{d}=\{l+1\}}\prod_{i=1}^{d}(|\beta_{i}|-1)!\varphi_{\sum_{i\in\beta_{1}}\eta_{i}}\cdots\varphi_{\sum_{i\in\beta_{d}}\eta_{i}}\varphi_{\eta_{l+1}}\label{term_one}\\
 &+\sum_{d=1}^{l}(-1)^{l+1-d}\sum_{v=1}^{d} \sum_{\beta^{v}\in\pi(l+1,d)}\prod_{i=1}^{d}(|\beta_{i}^{v}|-1)!\varphi_{\sum_{i\in\beta_{1}^{v}}\eta_{i}}\cdots\varphi_{\sum_{i\in\beta_{d}^{v}}\eta_{i}}.\label{term_two}
 \end{align}
 Let us separate the terms associated with $d=l+1$ from (\ref{term_one}) and separate the terms related  to $d=1$ from (\ref{term_two}); then we combine other terms in (\ref{term_one}) and (\ref{term_two}). Therefore, we have 
 \begin{align*}
 p_{(\eta_{1},\cdots,\eta_{l},\eta_{l+1})}^{o}
 =&\varphi_{\eta_{1}}\cdots\varphi_{\eta_{l}}\varphi_{\eta_{l+1}}
 +\sum_{d=2}^{l}(-1)^{l+1-d}\Bigg[\sum_{\beta\in\pi(l+1,d),\beta_{d}=\{l+1\}}\\
 &\prod_{i=1}^{d}(|\beta_{i}|-1)!\varphi_{\sum_{i\in\beta_{1}}\eta_{i}}\cdots\varphi_{\sum_{i\in\beta_{d}}\eta_{i}}\varphi_{\eta_{l+1}}\\
 &+\sum_{v=1}^{d} \sum_{\beta^{v}\in\pi(l+1,d)}(|\beta_{1}^{v}|-1)!\cdots(|\beta_{d}^{v}|-1)!\varphi_{\sum_{i\in\beta_{1}^{v}}\eta_{i}}\cdots\varphi_{\sum_{i\in\beta_{d}^{v}}\eta_{i}}\Bigg]\\
 &+(-1)^{l+1-1}(l+1-1)!\varphi_{\sum_{i=1}^{l+1}\eta_{i}}\\
=&\sum_{d=1}^{l+1}(-1)^{l+1-d}\sum_{\beta\in\pi(l+1,d)}\prod_{i=1}^{d}(|\beta_{i}|-1)!\varphi_{\sum_{i\in\beta_{1}}\eta_{i}}\cdots\varphi_{\sum_{i\in\beta_{d}}\eta_{i}}.
 \end{align*}
 The proof is thus completed.

 \end{proof}

\bibliography{final}{}

\def\cprime{$'$} \def\cprime{$'$} \def\cprime{$'$} \def\cprime{$'$}
  \def\cprime{$'$}
\providecommand{\bysame}{\leavevmode\hbox to3em{\hrulefill}\thinspace}
\providecommand{\MR}{\relax\ifhmode\unskip\space\fi MR }
\providecommand{\MRhref}[2]{%
  \href{http://www.ams.org/mathscinet-getitem?mr=#1}{#2}
}
\providecommand{\href}[2]{#2}
\begin{thebibliography}{10}

\bibitem{MR1649005}
Donald~A. Dawson and Shui Feng, \emph{Large deviations for the {F}leming-{V}iot
  process with neutral mutation and selection}, Stochastic Process. Appl.
  \textbf{77} (1998), no.~2, 207--232. \MR{1649005 (99h:60050)}

\bibitem{MR1815182}
\bysame, \emph{Large deviations for the {F}leming-{V}iot process with neutral
  mutation and selection. {II}}, Stochastic Process. Appl. \textbf{92} (2001),
  no.~1, 131--162. \MR{1815182 (2001m:60060)}

\bibitem{MR2571413}
Amir Dembo and Ofer Zeitouni, \emph{Large deviations techniques and
  applications}, Stochastic Modelling and Applied Probability, vol.~38,
  Springer-Verlag, Berlin, 2010, Corrected reprint of the second (1998)
  edition. \MR{2571413 (2011b:60094)}

\bibitem{MR1174426}
S.~N. Ethier, \emph{Eigenstructure of the infinitely-many-neutral-alleles
  diffusion model}, J. Appl. Probab. \textbf{29} (1992), no.~3, 487--498.
  \MR{1174426 (93m:60156)}

\bibitem{MR1235429}
S.~N. Ethier and R.~C. Griffiths, \emph{The transition function of a
  {F}leming-{V}iot process}, Ann. Probab. \textbf{21} (1993), no.~3,
  1571--1590. \MR{1235429 (95a:60101)}

\bibitem{MR615945}
S.~N. Ethier and Thomas~G. Kurtz, \emph{The infinitely-many-neutral-alleles
  diffusion model}, Adv. in Appl. Probab. \textbf{13} (1981), no.~3, 429--452.
  \MR{615945 (82j:60143)}

\bibitem{MR1205982}
\bysame, \emph{Fleming-{V}iot processes in population genetics}, SIAM J.
  Control Optim. \textbf{31} (1993), no.~2, 345--386. \MR{1205982 (94d:60131)}

\bibitem{MR1302692}
\bysame, \emph{Convergence to {F}leming-{V}iot processes in the weak atomic
  topology}, Stochastic Process. Appl. \textbf{54} (1994), no.~1, 1--27.
  \MR{1302692 (95m:60075)}

\bibitem{MR2358634}
Shui Feng, \emph{Large deviations associated with {P}oisson-{D}irichlet
  distribution and {E}wens sampling formula}, Ann. Appl. Probab. \textbf{17}
  (2007), no.~5-6, 1570--1595. \MR{2358634 (2009b:60088)}

\bibitem{MR2663265}
\bysame, \emph{The {P}oisson-{D}irichlet distribution and related topics},
  Probability and its Applications (New York), Springer, Heidelberg, 2010,
  Models and asymptotic behaviors. \MR{2663265 (2012c:60138)}

\bibitem{MR2737405}
Shui Feng, Wei Sun, Feng-Yu Wang, and Fang Xu, \emph{Functional inequalities
  for the two-parameter extension of the infinitely-many-neutral-alleles
  diffusion}, J. Funct. Anal. \textbf{260} (2011), no.~2, 399--413. \MR{2737405
  (2011k:60267)}

\bibitem{MR1887170}
Shui Feng and Jie Xiong, \emph{Large deviations and quasi-potential of a
  {F}leming-{V}iot process}, Electron. Comm. Probab. \textbf{7} (2002), 13--25
  (electronic). \MR{1887170 (2002m:60047)}

\bibitem{MR0526801}
J.~F.~C. Kingman, \emph{Random partitions in population genetics}, Proc. Roy.
  Soc. London Ser. A \textbf{361} (1978), no.~1704, 1--20. \MR{0526801 (58
  \#26167)}

\bibitem{MR509954}
\bysame, \emph{The representation of partition structures}, J. London Math.
  Soc. (2) \textbf{18} (1978), no.~2, 374--380. \MR{509954 (80a:05018)}

\bibitem{MR0368264}
J.~F.~C. Kingman, S.~J. Taylor, A.~G. Hawkes, A.~M. Walker, David~Roxbee Cox,
  A.~F.~M. Smith, B.~M. Hill, P.~J. Burville, and T.~Leonard, \emph{Random
  discrete distribution}, J. Roy. Statist. Soc. Ser. B \textbf{37} (1975),
  1--22, With a discussion by S. J. Taylor, A. G. Hawkes, A. M. Walker, D. R.
  Cox, A. F. M. Smith, B. M. Hill, P. J. Burville, T. Leonard and a reply by
  the author. \MR{0368264 (51 \#4505)}

\bibitem{MR1711591}
Zenghu Li, Tokuzo Shiga, and Lihua Yao, \emph{A reversibility problem for
  {F}leming-{V}iot processes}, Electron. Comm. Probab. \textbf{4} (1999),
  65--76 (electronic). \MR{1711591 (2001e:60097)}

\bibitem{MR0368829}
G.~A. Watterson, \emph{Models for the logarithmic species abundance
  distributions}, Theoret. Population Biology \textbf{6} (1974), 217--250.
  \MR{0368829 (51 \#5067)}

\bibitem{MR2184086}
Kai-Nan Xiang and Tu-Sheng Zhang, \emph{Small time asymptotics for
  {F}leming-{V}iot processes}, Infin. Dimens. Anal. Quantum Probab. Relat. Top.
  \textbf{8} (2005), no.~4, 605--630. \MR{2184086 (2007b:60065)}

\end{thebibliography}
\bibliographystyle{amsplain}
\end{document}